\documentclass[12pt,reqno]{amsart}

\usepackage[a4paper, left=3cm, right=3cm, bottom=2.5cm]{geometry}

\usepackage[english]{babel}
\usepackage[latin1]{inputenc}
\usepackage[T1]{fontenc}
\usepackage{amssymb,amsmath,amstext,amsfonts,amsthm,braket}
\usepackage{mathtools}
\usepackage{verbatim}
\usepackage{relsize}
\usepackage[colorlinks=true, urlcolor=blue, linkcolor=blue, citecolor=blue]{hyperref}
\usepackage{fancyvrb}
\usepackage{tikz-cd}
\usepackage[percent]{overpic}
\usepackage{comment}
\usepackage{algorithm}
\usepackage{algpseudocode}
\usepackage{todonotes}
\usepackage{xcolor}
\usepackage{cleveref}
\usepackage{gensymb}
\usepackage{graphicx}
\usepackage{bbold}
\usepackage{bm}
\usepackage[super]{nth}
\algrenewcommand{\algorithmiccomment}[1]{\hfill\textcolor{gray}{\texttt{\#} #1}}

\usepackage[normalem]{ulem}
\usepackage{enumitem}
\usepackage{neuralnetwork}
\usepackage{hyperref}
\usepackage{booktabs}

\numberwithin{equation}{section}
\theoremstyle{plain}
\newtheorem*{theorem*}{Theorem}
\newtheorem{theorem}{Theorem}
\numberwithin{theorem}{section}
\newtheorem{proposition}[theorem]{Proposition}
\newtheorem{lemma}[theorem]{Lemma}
\newtheorem{corollary}[theorem]{Corollary}
\newtheorem{conjecture}[theorem]{Conjecture}

\theoremstyle{definition}
\newtheorem{definition}[theorem]{Definition}
\newtheorem{remark}[theorem]{Remark}
\newtheorem{example}[theorem]{Example}

\newcommand{\C}{\mathbb{C}}

\newcommand{\Z}{\mathbb{Z}}
\newcommand{\R}{\mathbb{R}}
\newcommand{\N}{\mathbb{N}}
\renewcommand{\P}{\mathbb{P}}
\newcommand{\V}{\mathcal{V}}
\newcommand{\M}{\mathcal{M}}
\renewcommand{\O}{\mathcal{O}}
\renewcommand{\SS}{\mathcal{S}}
\newcommand{\QQ}{\mathcal{Q}}
\newcommand{\T}{\mathcal{T}}

\newcommand{\x}{\mathbf{x}}
\newcommand{\y}{\mathbf{y}}
\newcommand{\w}{\mathbf{w}}
\newcommand{\ttheta}{\bm{\theta}}
\renewcommand{\d}{\mathbf{d}}

\newcommand{\del}{\partial}
\DeclareMathOperator{\Sym}{Sym}
\DeclareMathOperator{\gr}{Gr}
\DeclareMathOperator{\im}{im}

\DeclareMathOperator{\Gr}{Gr}
\DeclareMathOperator{\edim}{edim}
\DeclareMathOperator{\ldeg}{Ldeg}

\definecolor{taskcol}{rgb}{0.8, 0.5, 0.95}
\definecolor{notecol}{rgb}{0.95, 0.8, 0.5}

\DeclareMathOperator{\rk}{rk}

\title{Geometry of polynomial neural networks}

\author{Kaie Kubjas}
\address{Department of Mathematics and Systems Analysis, Aalto University, Espoo}
\email{kaie.kubjas@aalto.fi}

\author{Jiayi Li}
\address{Department of Statistics \& Data Science, University of California, Los Angeles}
\email{jiayi.li@g.ucla.edu}

\author{Maximilian Wiesmann}
\address{Max Planck Institute for Mathematics in the Sciences, Leipzig}
\email{wiesmann@mis.mpg.de}

\begin{document}
\maketitle

\begin{abstract}

    We study the expressivity and learning process for polynomial neural networks (PNNs) with monomial activation functions. The weights of the network parametrize the neuromanifold. In this paper,\ we study certain neuromanifolds using tools from algebraic geometry: we give explicit descriptions as semialgebraic sets and characterize their Zariski closures, called neurovarieties. We study their dimension and associate an algebraic degree, the learning degree, to the neurovariety. The dimension serves as a geometric measure for the expressivity of the network, the learning degree is a measure for the complexity of training the network and provides upper bounds on the number of learnable functions. These theoretical results are accompanied with experiments. 
    
\end{abstract}



\section{Introduction}

Over the past decade, neural networks have achieved remarkable success, primarily driven by advancements in deep learning. With the increased computational power, availability of large datasets, and algorithmic innovations, deep learning has led to groundbreaking achievements in areas like image and speech recognition, natural language processing, and autonomous systems.\ In addition, deep neural networks outperform many traditional statistical models, significantly impacting fields such as healthcare and finance.\ Despite the empirical triumph, the underlying theoretical behaviors of neural networks remain an open and active field of research. 

In particular, researchers study various activation functions to understand their role in introducing non-linearity, affecting gradient propagation, and influencing computational efficiency.\ Polynomial activation functions have gained interest for their ability to introduce higher-order interactions between inputs, allowing networks to model complex, non-linear phenomena more effectively \cite{oh2003polynomial}.  Although feedforward neural networks with polynomial activation functions are well-known to be non-universal approximators \cite{hornik1989multilayer},  there are many practical tasks where architectures with polynomial activations have outperformed other ones, especially in environments where data relationships are polynomial in nature.\ In particular, polynomial neural networks have led to state-of-the-art results in engineering tasks like face detection from cluttered images \cite{huang2003face}, image generation \cite{chrysos2020p},   3D shape recognition \cite{yavartanoo2021polynet}, as well as  financial applications like forecasting trading signals \cite{ghazali2011dynamic},  uncertain natural frequency quantification \cite{dey2016uncertain}, and estimating stock closing indices \cite{nayak2018estimating}.\ However, the choice of degree and the potential for overfitting are critical considerations in applying polynomial neural networks to individual practical tasks.

Compared to other low-degree activation functions, polynomial functions capture intricate patterns within data without the need for additional layers, potentially reducing model complexity and computational costs.\ In addition, common neural network activation functions, including sigmoid and ReLU, can be effectively approximated using polynomial ratios.\ Recent work also shows that fully-connected feedforward neural networks using ratios of polynomials as activation functions approximate smooth functions more efficiently than ReLU networks \cite{boulle2020rational}.\ Our exploration in the realm of polynomial networks lays the groundwork for further investigation of other activation functions.

In this paper, we perform an algebro-geometric study of neuromanifolds\footnote{By an abuse of terminology we use the term \emph{neuromanifold} to denote the functional space of a neural network, although this space is almost never a smooth manifold.} and their Zariski closures, neurovarieties, following the previous work by \cite{kileel2019expressive}.\ 
In~\S \ref{sec:background}, we give the background on polynomial neural networks, neuromanifolds and neurovarieties. Neuromanifolds provide a natural generalization of the set of symmetric tensors of bounded symmetric rank.  In~\S \ref{sec:neuromanifolds}, we recall the connection between homogeneous polynomials and symmetric tensor decompositions, and characterize neuromanifolds for some shallow polynomial neural networks. In~\S \ref{sec:neurovarities}, we describe different approaches for studying neurovarieties. The dimension of the neurovariety provides a measure for the expressivity of the network. In~\S \ref{sec:dimension}, we study these dimensions. For a shallow polynomial neural network with a single output unit, this corresponds to the Alexander--Hirschowitz Theorem in neural network settings. In the deep case, we present conjectures supported by empirical evidence. Finally, in~\S \ref{sec:optimization}, we study the optimization process of a polynomial neural network: from a static perspective, we describe the complexity of the optimization landscape by its \emph{learning degree} and compute it for a family of architectures; for the dynamic optimization process, we review the backpropagation algorithm and summarize known results on linear neural networks from previous literature in both the machine learning and algebraic statistics communities.\par 
We provide code for computations and experiments at \cite{mathrepo}.

\textbf{Main contributions.}  For those new to the topic, we recommend coming back to this paragraph on the main contributions after reading~\Cref{sec:background}.  We characterize the neuromanifold $\M_{(d_0,d_1,1),2}$ in~\Cref{lemma:d0d11r2}, $\M_{(d_0,1,d_2),2}$ in~\Cref{cor:dd01d2_r2},  $\M_{(d_0,1,d_2),r}$ in~\Cref{lemma:dd01d2_r}, and 
 $\M_{(d_0, d_1, d_2), r}$ with $d_1 \geq \binom{r+d_0-1}{r}$ in~\Cref{prop:fillingMfd}. We characterize the neurovariety $\V_{(2,2,d_2),2}$ in~\Cref{prop:d22k_r2} and partially the neuromanifold $\M_{(2,2,d_2),2}$ in~\Cref{prop:semialgebraic_description_d22k_r2} by studying the fibers of the PNN, the neurovariety $\V_{(3,3,3),2}$ in~\Cref{example:grassm} by using Grassmannians and the neurovariety $\V_{(2,2,2,2),2}$ in~\Cref{example:Hilbert-Burch} by applying the Hilbert--Burch Theorem. In~\Cref{prop:neuromanifold-with-width-1-layer}, we show that if a PNN has a layer of width one, then the widths of layers after the width-$1$ layer except for the output layer do not influence the neuromanifold. In~\Cref{cor:defective-subnetwork-implies-defective-full-network}, we show that if a PNN has a defective subnetwork, then it is itself defective unless it is filling. In~\Cref{sec:symmetries-of-an-exceptional-shallow-network}, we describe a family of non-trivial symmetries for the neuromanifold $\V_{(5,7,1),3}$ that does not have expected dimension. We introduce the learning degree, a characteristic of the optimization landscape, in Definition \ref{def:learningDeg}. In~\Cref{thm:learning-degree-is-upper-bounded-by-the-generic-ED-degree}, we show that learning degree of a PNN with respect to the Euclidean loss function exists and it is bounded above by the generic Euclidean distance degree of its neurovariety. In~\Cref{prop:EDded22k}, we prove that the learning degree of $\V_{(2,2,k),2}$ for $k\geq 2$ is at most $8k^2 - 12k + 3$.

\textbf{Acknowledgments.} This paper was written in connection with the ``Apprenticeship Week: Varieties from Statistics'' which took place at IMSI as part of the long program ``Algebraic Statistics and Our Changing World''. We would like to thank Miles Bakenhus and Maksym Zubkov who participated in the work of this group part of the time. \Cref{lem:nn2} is due to Maksym Zubkov. We would like to thank Elizabeth Gross, Leonie Kayser, Joe Kileel, Mark Kong, Guido Mont{\'u}far, Alessandro Oneto, Bernd Sturmfels, and Matthew Trager for useful discussions and various suggestions. We thank Lisa Seccia and Nils Sturma for careful reviews of an earlier version of this paper. We also thank Giovanni Luca Marchetti and Vahid Shahverdi for pointing out a mistake in Section \ref{sec:learning-degree} in a previous version of this article. Part of this research was performed while the authors were visiting the Institute for Mathematical and Statistical Innovation (IMSI), which is supported by the National Science Foundation (Grant No.\ DMS-1929348). 

\section{Background} \label{sec:background}

A general $L$-layer \emph{(feedforward) neural network} is a composition of $L$ affine-linear maps with coordinatewise nonlinearity in between \cite{haykin1998neural}. Let $F_{\ttheta}$ be a feedforward neural network with  parameters $\ttheta$, 
\[F_{\ttheta}(\x)=f_L\circ\sigma_{L-1}\circ f_{L-1}\circ\cdots \circ f_2\circ\sigma_{1}\circ f_1(\x),\]
where \[f_{l}(\x)\colon \R^{d_{l-1}}\rightarrow \R^{d_l},\quad \x\mapsto W_l\x+\textbf{b}_l\]
and the function $\sigma_l \colon \R^{d_l}\rightarrow \R^{d_l}$ is the \emph{activation function}, typically a real function that is applied coordinatewise, i.e.,\ $\sigma_l = (\sigma_{l,1},\dots, \sigma_{l,d_l})$ with $\sigma_{l,j}\colon\R\rightarrow\R$. The matrices $W_1,\ldots,W_L$ are the \emph{weights} of the neural network, and $b_1,\ldots,b_L$ are referred to as the \emph{biases}. Together, they constitute the parameter set $\ttheta$. \par

\begin{figure}
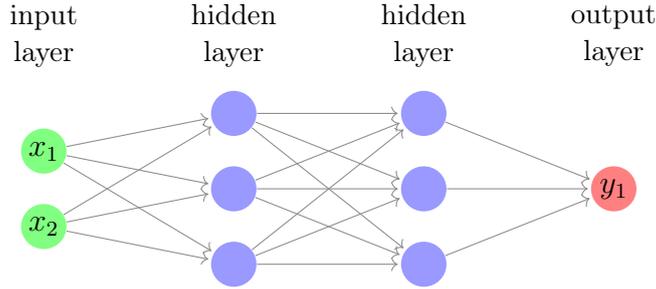

    \begin{neuralnetwork}[height=3]
        \newcommand{\xx}[2]{${x}_#2$}
        \newcommand{\yy}[2]{${y}_#2$}
        \inputlayer[count=2, bias=false, title={\small input\\layer}, text=\xx]
        \hiddenlayer[count=3, bias=false, title={\small hidden\\layer}] \linklayers
        \hiddenlayer[count=3, bias=false, title={\small hidden\\layer}] \linklayers
        \outputlayer[count=1, title={\small output\\layer}, text=\yy] \linklayers
    \end{neuralnetwork}
    \caption{A neural network architecture with widths $\d = (2,3,3,1)$, input $\x = (x_1, x_2)^T$ and output $y_1$.}
    \label{fig:NNarch}
\end{figure}

Pictorially, a neural network architecture can be represented as in Figure \ref{fig:NNarch}. Each node is called a \emph{neuron} and each column of neurons forms a \emph{layer} of the network. The first layer is called \emph{input layer}, the last layer is called \emph{output layer} and the remaining layers are referred to as \emph{hidden layers}. A neural network with exactly one hidden layer is called \emph{shallow}. The number of neurons in the $i^{\text{th}}$ layer is the $i^{\text{th}}$ \emph{width}, denoted $d_i$. The vector of widths $\d =(d_0,d_1,\dots,d_L)$ together with a choice of activation functions $\bm{\sigma} = (\sigma_1,\sigma_2,\dots, \sigma_{L-1})$ constitute the \emph{architecture} of the network.Arrows between layers represent the maps $f_{l}$. Note that all architectures in this paper are considered to be fully-connected feedforward neural networks.  \par 
The goal of deep learning is to approximate a target function $f\colon\R^{d_0}\rightarrow\R^{d_L}$ with a neural network $F_{\ttheta}$ of a chosen architecture $(\d = (d_0,d_1,\dots,d_L),\bm{\sigma})$. This amounts to an optimization task over the space of parameters, the ``learning'' or ``training'' process. For more details, see \S \ref{sec:optimization}.\par 
The \emph{parameter map}
\[
    \Psi_{\d,\bm{\sigma}}\colon \R^N \rightarrow \mathrm{Fun}(\R^{d_0}, \R^{d_L}),\quad \ttheta\mapsto F_{\ttheta}
\]
associates a tuple of parameters $\ttheta$ with the corresponding neural network $F_{\ttheta}$. Its image is called the \emph{neuromanifold} $\M_{\d,\bm{\sigma}}$ and consists of all functions a network with architecture $(\d,\bm{\sigma})$ can learn. The neuromanifold is also referred to as \emph{functional space} in the literature. Note that oftentimes $\M_{\d,\bm{\sigma}}$ is not a smooth manifold.\par 
In practice, typical choices of activation functions are \emph{ReLU} or \emph{sigmoid} functions
\[
    \sigma_{\text{ReLU}}(x) = \max\{0,x\},\quad \sigma_{\text{sigmoid}}(x) = \frac{e^x}{1+e^x}.
\]
In this paper, we focus on \emph{monomial} activations $x\mapsto x^r$. This turns the network $F_{\ttheta}$ into a polynomial map and makes the study of these networks amenable to techniques from algebraic geometry. \par
The map $f_l$ can be written as $f_l(\x) = W_{l}\x + \mathbf{b}_l$, where $W_l\in \R^{d_l\times d_{l-1}}$ is a linear map and $\mathbf{b}_l\in \R^{d_i}$. We will be considering networks without biases to make the polynomial map $F_{\ttheta}$ homogeneous. Let us summarize the setup in the following.

\begin{definition}
    A \emph{polynomial neural network (PNN)} $p_{\w}$ with architecture $(\d=(d_0,d_1,\dots,d_L),r)$ is a feedforward neural network
    \[
        p_{\w} = W_L\circ \sigma_{L-1}\circ W_{L-1} \circ \sigma_{L-2} \circ \dots \circ \sigma_1\circ W_1 \colon \R^{d_0}\rightarrow\R^{d_L}
    \]
    where $W_i\in\R^{d_i\times d_{i-1}}$ are linear maps and the activation functions 
    \[
        \sigma_i(\x) = \rho_r(\x) \coloneqq (x_1^r,\dots, x_n^r)
    \] 
    are monomial. The number $r$ is called the \emph{activation degree} of $p_{\w}$. The parameters $\w$ are given by the entries of the matrices $W_i$, i.e.,\ $\w = (W_1,W_2,\dots,W_L)$.
\end{definition}

A PNN $p_{\w}$ with architecture $(\d,r)$ is a homogeneous polynomial map of degree $r^{L-1}$. Hence, the associated parameter map $\Psi_{\d,r}$ maps an $L$-tuple of matrices $(W_1,W_2,\dots,W_L)$ to a $d_L$-tuple of homogeneous polynomials of degree $r^{L-1}$ in $d_0$ variables, i.e.,\
\[
    \Psi_{\d, r}: \R^{d_1 \times d_0} \times \cdots \times \R^{d_L \times d_{L-1}} \rightarrow (\Sym_{r^{L-1}}(\R^{d_0}))^{d_L}, \quad \w \mapsto p_{\w}=
    \begin{pmatrix}
        p_{\w}^{(1)}\\
        \vdots\\
        p_{\w}^{(d_L)}
    \end{pmatrix}.
\]

We can identify elements in the image with their vectors of coefficients in 
\[
    \R^{d_L\binom{r^{L-1}+d_0-1}{d_0-1}}.
\]
Coordinates of this space will be denoted by $c^{(j)}_I$ so that $c^{(j)}_I$ is the coefficient of the monomial $\x^I$ in the polynomial $p_{\w}^{(j)}(\x)$; here, $I$ is a multiindex in 
\[
    \binom{[d_0]}{r^{L-1}} = \left\{r^{L-1}\text{ element subsets in }\{1,2,\dots,d_0\}\right\}.
\]

\begin{definition}
    The image of $\Psi_{\d,r}$ is the \emph{neuromanifold} $\M_{\mathbf{d},r}$. Its Zariski closure is the \emph{neurovariety} $\V_{\mathbf{d},r}$, an affine variety inside $(\Sym_{r^{L-1}}(\R^{d_0}))^{d_L}$.
\end{definition}

The neuromanifold $\Psi_{\d,r}$ is a semialgebraic set inside $(\Sym_{r^{L-1}}(\R^{d_0}))^{d_L}$, which means that it is a finite union of sets that can be defined by polynomial equalities and inequalities. The observation that the neuromanifold $\Psi_{\d,r}$ is a semialgebraic set follows from the Tarski--Seidenberg theorem which states that a polynomial image of a semialgebraic set is again semialgebraic.

Increasing the widths of hidden layers gives a containment of neuromanifolds as shown in the following proposition.

\begin{proposition}
    Let $\mathbf{d}=(d_0,\dots,d_i,\dots,d_L)$ and let $\d'=(d_0,\dots,d'_i,\dots,d_L)$ be a tuple which differs from $\d$ precisely in the $i^{\text{th}}$ entry for $0<i<L$ and assume $d'_i \geq d_i$. Then $\M_{\mathbf{d},r}\subseteq \M_{\d',r}$, in particular $\V_{\mathbf{d},r}\subseteq \V_{\d',r}$.
\end{proposition}

\begin{proof}
Let $\w = (W_1,\ldots,W_L)$ be a parameter vector for the architecture $\d$ and $p_{\w} \in \M_{\mathbf{d},r}$ the corresponding polynomial network. Let $\w' = (W'_1,\ldots,W'_L)$ be the parameter vector for the architecture $\d'$ such that each $W'_i$ has $W_i$ as left-top submatrix and zeros elsewhere. Then $p_{\w}=p_{\w'} \in \M_{\d',r}$.
\end{proof}

As $\M_{\d,r}$ is semialgebraic, it has the same dimension as its Zariski closure, $\V_{\d,r}$. In \cite{kileel2019expressive}, the dimension of $\V_{\d,r}$ was proposed as a measure for the \emph{expressivity} of the network architecture $(\d,r)$.

\begin{definition}
    An architecture $(\d,r)$ is \emph{filling} if $\V_{\mathbf{d},r} = (\Sym_{r^{L-1}}(\R^{d_0}))^{d_L}$. In this case, we say that $\M_{\mathbf{d},r}$ is \emph{thick}, i.e., it has positive Lebesgue measure.
\end{definition}

The case of filling architectures is particularly interesting from a machine learning perspective as filling networks have the most expressive power: if $\M_{\d,r}=(\Sym_{r^{L-1}}(\R^{d_0}))^{d_L}$, any target function in $(\Sym_{r^{L-1}}(\R^{d_0}))^{d_L}$ can potentially be learned \emph{exactly} by the network. In the case of non-filling architectures, a general target function can only be approximated by the network. For more details on the learning process, see \S\ref{sec:optimization}.  Moreover, from an optimization perspective it is advantageous to work with thick neuromanifolds as non-thick neuromanifolds are known to be non-convex, see \cite[Proposition 7]{kileel2019expressive}.

\begin{example}
For architecture $\d = (2,1,1),\ r = 2$ with input $\x=(x_1,\ x_2)^T$ and parameters
		\[W_1 = \begin{pmatrix}
			w_{111} & w_{112} 
		\end{pmatrix},\ W_2 = \begin{pmatrix} 
		w_{211}
	\end{pmatrix},\]
		the network is
		\[p_{\w}(\x) =W_2\rho_2W_1 \x =\begin{pmatrix}
			w_{211}(w_{111}x_1+w_{112}x_2)^2
		\end{pmatrix}.\]
    \begin{figure} 
    \centering
    \includegraphics[width=0.4\textwidth]{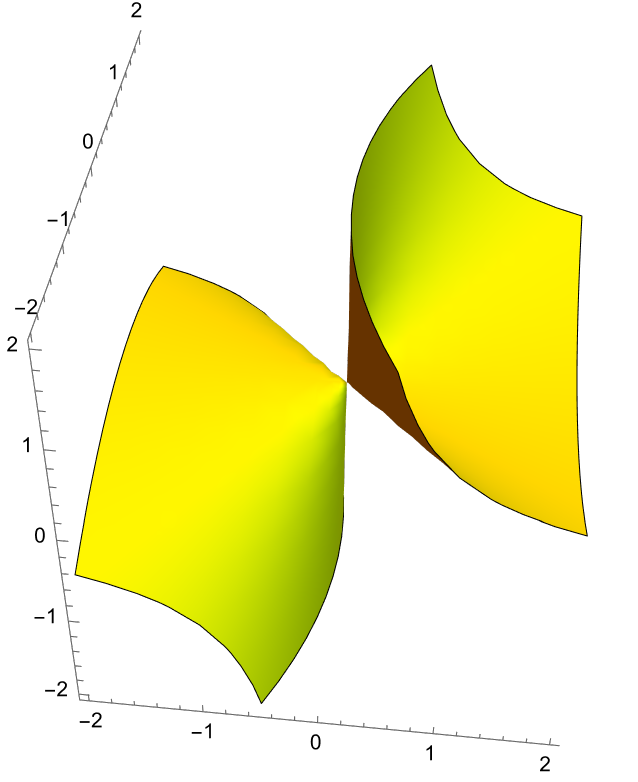}
    \caption{The neuromanifold $\M_{(2,1,1), 2}$ in $\text{Sym}_2(\R^2)\cong \R^3$.}
    \label{fig:d211_r2}
\end{figure}
		There are three parameters in $W_1$ and $W_2$ which $\Psi_{(2,1,1),2}$ maps to the quadric above with $\dim(\Sym_2(\R^2)) = 3$ coefficients. Let $c_{11}, c_{12}$ and $c_{22}$ be coordinates for this space, representing the coefficients of $x_1^2, x_1x_2$ and $x_2^2$, respectively. The neuromanifold is a hypersurface defined by the quadratic equation $c_{12}^2-c_{11}c_{22}=0$, which is shown in~\Cref{fig:d211_r2}. The neurovariety is equal to the neuromanifold, and hence it is not filling. See Lemma \ref{lemma:d0d11r2} for more details and generalizations.
\end{example}

\begin{example}[\cite{kileel2019expressive}, Example 3]
		For architecture $\d = (2,2,3),\ r = 2$ with input $\x=(x_1,\ x_2)^T$ and parameters
		\[W_1 = \begin{pmatrix}
			w_{111} & w_{112} \\
			w_{121} & w_{122} 
		\end{pmatrix}, \ W_2 = \begin{pmatrix}
		w_{211} & w_{212} \\ 
		w_{221} & w_{222} \\
		w_{231} & w_{232} 
	\end{pmatrix},\]
		the network is
		\[p_{\w}(\x) =W_2\rho_2W_1 \x =\begin{pmatrix}
			w_{211}(w_{111}x_1+w_{112}x_2)^2 + w_{212}(w_{121}x_1+w_{122}x_2)^2 \\
			w_{221}(w_{111}x_1+w_{112}x_2)^2 + w_{222}(w_{121}x_1+w_{122}x_2)^2 \\
			w_{231}(w_{111}x_1+w_{112}x_2)^2 + w_{232}(w_{121}x_1+w_{122}x_2)^2 
		\end{pmatrix}.\]
		There are ten parameters in $W_1$ and $W_2$ which $\Psi_{(2,2,3),2}$ maps to the triple of quadrics above with $\dim(\Sym_2(\R^2)^3) = 9$ coefficients. The neurovariety is a hypersurface defined by the cubic equation
		\[\det \begin{pmatrix}
			c_{11}^{(1)} & c_{12}^{(1)}  & c_{22}^{(1)} \\
			c_{11}^{(2)} & c_{12}^{(2)}  & c_{22}^{(2)} \\
			c_{11}^{(3)} & c_{12}^{(3)}  & c_{22}^{(3)} 
		\end{pmatrix}=0.\]
		 This variety has dimension eight, implying the architecture $\d = (2,2,3), r=2$ is not filling. For more details, see Proposition \ref{prop:d22k_r2} below. Note that in this case $\M_{\d,r}\subsetneq \V_{\d,r}$, see Proposition \ref{prop:semialgebraic_description_d22k_r2}.
\end{example}

There is a na\"ive expectation for the dimension of the neurovariety, namely the number of parameters. However, one immediately observes a symmetry in the parameters for any network architecture, called \emph{multi-homogeneity}.

\begin{lemma}[{\cite[Lemma 13]{kileel2019expressive}}]
    \label{lem:multiHomog}
    For all diagonal matrices $D_i\in\R^{d_i\times d_i}$ and permutation matrices $P_i\in\Z^{d_i\times d_i}$ ($i=1,\dots,L-1$), the parameter map $\Psi_{\d,r}$ returns the same network under the replacement:
    \[
        \begin{array}{l}
            W_1 \gets P_1D_1W_1 \\
            W_2 \gets P_2D_2W_2D_1^{-r}P_1^T \\
            W_3 \gets P_3D_3W_3D_2^{-r}P_2^T \\
            \hspace{2em}\vdots \\
            W_L \gets W_LD_{L-1}^{-r}P_{L-1}^T.
        \end{array}
    \]
\end{lemma}
Hence, a generic fiber of $\Psi_{\d,r}$ has dimension at least $\sum_{i=1}^{L-1}d_i$. We call the number of parameters subtracted by this dimension the expected dimension of the neurovariety.

\begin{definition}\label{def:expected-dim}
    We define the \emph{expected dimension} of the neurovariety $\V_{\d,r}$ to be
    \[
        \edim(\V_{\d,r}) \coloneqq \min\left\{d_L + \sum_{i=0}^{L-1}(d_id_{i+1} - d_{i+1}), d_L\binom{d_0 + r^{L-1} - 1}{r^{L-1}}\right\}.
    \]
    By Lemma \ref{lem:multiHomog}, $\dim(\V_{\d,r}) \leq \edim(\V_{\d,r})$. The difference $\edim(\V_{\d,r}) - \dim(\V_{\d,r})$ is the \emph{defect} of $\V_{\d,r}$. If the defect is nonzero $\V_{\d,r}$ is called \emph{defective}. We refer to 
    \[
        \dim((\Sym_{r^{L-1}}(\R^{d_0}))^{d_L}) = d_L\binom{d_0 + r^{L-1} - 1}{r^{L-1}}
    \]
    as the \emph{ambient dimension} of $\V_{\d,r}$.
\end{definition}

In~\Cref{table:small_examples} we compute the ideal of the neurovariety and its dimension for shallow polynomial neural networks with $d_i \in [3]$ and $r=2$. We also compare the neurovariety with the neuromanifold.

\begin{table}
\begin{center}
\footnotesize
\begin{tabular}{ c c c c c c}
 $\d$ & dim & edim & amb dim & ideal & $\M_{\d,2} = \V_{\d,2}$?\\ \hline 
(1,1,1) & 1 & 1 & 1 & $\langle 0\rangle$ & yes \\
(1,1,2) & 2 & 2 & 2 & $\langle 0\rangle$ & yes \\
(1,1,3) & 3 & 3 & 3 & $\langle 0\rangle$ & yes \\ \\
(1,2,1) & 1 & 1 & 1 & $\langle 0\rangle$ & yes \\ 
(1,2,2) & 2 & 2 & 2 & $\langle 0\rangle$ & yes \\
(1,2,3) & 3 & 3 & 3 & $\langle 0\rangle$ & yes \\ \\

(1,3,1) & 1 & 1 & 1 & $\langle 0\rangle$ & yes \\
(1,3,2) & 2 & 2 & 2 & $\langle 0 \rangle$ & yes \\
(1,3,3) & 3 & 3 & 3 & $\langle 0 \rangle$ & yes \\ \\

(2,1,1) & 2 & 2 & 3 & determinantal,~\Cref{cor:dd01d2_r2}  & yes,~\Cref{cor:dd01d2_r2}\\
(2,1,2) & 3 & 3 & 6 & determinantal,~\Cref{cor:dd01d2_r2}  & yes,~\Cref{cor:dd01d2_r2} \\
(2,1,3) & 4 & 4 & 9 & determinantal,~\Cref{cor:dd01d2_r2}  & yes,~\Cref{cor:dd01d2_r2} \\ \\

(2,2,1) & 3 & 3 & 3 & $\langle 0 \rangle$ & yes, \Cref{lemma:d0d11r2}\\
(2,2,2) & 6 & 6 & 6 & $\langle 0 \rangle$ & no, \Cref{prop:semialgebraic_description_d22k_r2}\\
(2,2,3) & 8 & 8 & 9 & determinantal, Proposition~\ref{prop:d22k_r2} & no, \Cref{prop:semialgebraic_description_d22k_r2}\\ \\
 
(2,3,1) & 3 & 3 & 3 & $\langle 0 \rangle$ & yes, \Cref{lemma:d0d11r2} \\
(2,3,2) & 6 & 6 & 6 & $\langle 0 \rangle$ & yes, Proposition \ref{prop:fillingMfd} \\
(2,3,3) & 9 & 9 & 9 & $\langle 0 \rangle$ & yes, Proposition \ref{prop:fillingMfd}\\ \\

(3,1,1) & 3 & 3 & 6 &  determinantal,~\Cref{cor:dd01d2_r2} & yes,~\Cref{cor:dd01d2_r2}\\
(3,1,2) & 4 & 4 & 12 & determinantal,~\Cref{cor:dd01d2_r2} & yes,~\Cref{cor:dd01d2_r2}  \\
(3,1,3) & 5 & 5 & 18 & determinantal,~\Cref{cor:dd01d2_r2} & yes,~\Cref{cor:dd01d2_r2} \\ \\

(3,2,1) & 5  & 6 & 6 &  determinantal, \Cref{lemma:d0d11r2} & yes, \Cref{lemma:d0d11r2}\\
(3,2,2) & 8 & 8 & 12 & determinantal,  Example~\ref{example:table_unexplained_cases} & {no, Remark~\ref{rem:(d_0,2,d_2),r=2}} \\
(3,2,3) & 10 & 10 & 18 &  Example~\ref{example:table_unexplained_cases} & {no, Remark~\ref{rem:(d_0,2,d_2),r=2}} \\ \\

(3,3,1) & 6  &  6 & 6 & $\langle 0 \rangle$ & yes,~\Cref{cor:dd01d2_r2} \\
(3,3,2) & 12 & 12 & 12 & $\langle 0 \rangle$ & no, Lemma \ref{lem:nn2} \\
(3,3,3) & 15 & 15 & 18 &   \Cref{example:grassm}  & no, Lemma \ref{lem:nn2}
\end{tabular}
\end{center}
\caption{Properties of neurovarieties for shallow polynomial neural networks with $d_i \in [3]$ and $r=2$.}
\label{table:small_examples}
\end{table}

Most ideals of neurovarieties in~\Cref{table:small_examples} are described by results in the upcoming sections. The ideals for all but two cases ($\mathbf{d}=(3,2,3)$ and $\mathbf{d}=(3,3,3)$) in the table are determinantal. In general, we expect most neurovarieties to be non-determinantal. Among the table, only the neurovariety for the architecture $\mathbf{d}=(3,2,1)$ does not have expected dimension. The following example describes the cases that are not explained by any results in the rest of the paper.

\begin{example} \label{example:table_unexplained_cases}
The ideal of the neurovariety for the architecture $\mathbf{d}=(3,2,2),r=2$ is determinantal and it is generated by the $3 \times 3$ minors of
\[
\begin{pmatrix}
c_{11}^{(1)} & \frac{1}{2} c_{12}^{(1)} & \frac{1}{2} c_{13}^{(1)} & c_{11}^{(2)} & \frac{1}{2} c_{12}^{(2)} & \frac{1}{2} c_{13}^{(2)} \\
\frac{1}{2} c_{12}^{(1)} &  c_{22}^{(1)} & \frac{1}{2} c_{23}^{(1)} & \frac{1}{2} c_{12}^{(2)} &  c_{22}^{(2)} & \frac{1}{2} c_{23}^{(2)} \\
\frac{1}{2} c_{13}^{(1)} & \frac{1}{2} c_{23}^{(1)} &  c_{33}^{(1)}&  \frac{1}{2} c_{13}^{(2)} & \frac{1}{2} c_{23}^{(2)} & c_{33}^{(2)} 
\end{pmatrix}.
\]

The ideal of the neurovariety for the architecture $\d = (3,2,3),r=2$ contains the ideal generated by the $3 \times 3$ minors of the following $3 \times 9$ matrix:
\[
\begin{pmatrix}
c_{11}^{(1)} & \frac{1}{2} c_{12}^{(1)} & \frac{1}{2} c_{13}^{(1)} & c_{11}^{(2)} & \frac{1}{2} c_{12}^{(2)} & \frac{1}{2} c_{13}^{(2)} & c_{11}^{(3)} & \frac{1}{2} c_{12}^{(3)} & \frac{1}{2} c_{13}^{(3)} \\
\frac{1}{2} c_{12}^{(1)} &  c_{22}^{(1)} & \frac{1}{2} c_{23}^{(1)} & \frac{1}{2} c_{12}^{(2)} &  c_{22}^{(2)} & \frac{1}{2} c_{23}^{(2)} & \frac{1}{2} c_{12}^{(3)} &  c_{22}^{(3)} & \frac{1}{2} c_{23}^{(3)} \\
\frac{1}{2} c_{13}^{(1)} & \frac{1}{2} c_{23}^{(1)} &  c_{33}^{(1)}&  \frac{1}{2} c_{13}^{(2)} & \frac{1}{2} c_{23}^{(2)} & c_{33}^{(2)} & \frac{1}{2} c_{13}^{(3)} & \frac{1}{2} c_{23}^{(3)} & c_{33}^{(3)} 
\end{pmatrix}.
\]
However, it is not equal to this determinantal ideal. It is minimally generated by $94$ polynomials. The dimension of the determinantal variety is eleven while the neurovariety has dimension ten. 
\end{example}

\section{Neuromanifolds and Symmetric Tensor Decompositions} \label{sec:neuromanifolds}

In this section, we investigate the relationship between shallow PNNs and symmetric tensor decompositions. Building on this connection, we derive results for neuromanifolds for some shallow PNNs. This section has three subsections: \S \ref{sec:symmetric_tensor_decompositions} on symmetric tensors and homogeneous polynomials, \S \ref{sec:neuromanifolds_r2} on results for neuromanifolds for $r=2$ and \S \ref{sec:neuromanifolds_r>2} on results for neuromanifolds for general $r \in \N$. An overview of the results in \S \ref{sec:neuromanifolds_r2} and \S \ref{sec:neuromanifolds_r>2} can be found in~\Cref{table:neuromanifold_results}.

\begin{table}
\begin{center}
\footnotesize
\begin{tabular}{ c c c c}
$\d$ & $r$ & Assumptions & Result \\
\hline
$(d_0,d_1,1)$ & 2 & & \Cref{lemma:d0d11r2}\\
$(d_0,1,d_2)$ & 2 & & \Cref{cor:dd01d2_r2}\\
$(d_0,d_0,d_2)$ & 2 & $d_0,d_2 \geq 2$ & \Cref{lem:nn2}\\
$(d_0,1,d_2)$ & $r \in \N$ & & \Cref{lemma:dd01d2_r}\\
$(d_0, d_1, d_2)$ & $r \in \N$ & $d_1 \geq \binom{r+d_0-1}{r}$ & \Cref{prop:fillingMfd}\\
$(2,d_1,1)$ & $3$ & $d_1 \geq 2$ & \Cref{cor:neuromanifolds_from_tensors1}\\
$(2,d_1,1)$ & $4$ & $d_1 \geq 3$ & \Cref{cor:neuromanifolds_from_tensors1}\\
$(2,d_1,1)$ & $5$ & $d_1 \geq 3$ & \Cref{cor:neuromanifolds_from_tensors1}\\
$(3,d_1,1)$ & $4$ & $d_1\geq 6$ & \Cref{cor:neuromanifolds_from_tensors2}\\
$(3,d_1,1)$ & $5$ & $d_1\geq 7$ & \Cref{cor:neuromanifolds_from_tensors2}\\
$(4,d_1,1)$ & $3$ & $d_1\geq 5$ & \Cref{cor:neuromanifolds_from_tensors2}\\

\end{tabular}
\end{center}
\caption{An overview of architectures for which we give results about their neuromanifolds in \S \ref{sec:neuromanifolds}.}
\label{table:neuromanifold_results}
\end{table}

Consider a PNN with architecture $\mathbf{d}=(d_0,d_1,1)$. Then $\M_{\mathbf{d},r}$ consists of homogeneous polynomials of degree $r$ in $d_0$ many variables that can be represented as a linear combination of $d_1$ many linear forms raised to the $r^{\text{th}}$ power. We should emphasize here that the linear combinations are taken over the \emph{real} numbers. There is a bijection between the set of homogeneous polynomials of degree $r$ in $d_0$ many variables and the set of order-$r$ symmetric tensors of format $d_0 \times d_0 \times \cdots \times d_0$ as we will explain below. Formulated in the language of tensors, the neuromanifold $\M_{\mathbf{d},r}$ consists of order-$r$ symmetric tensors of format $d_0 \times d_0 \times \cdots \times d_0$ with \emph{real symmetric rank} $\leq d_1$.

\subsection{Symmetric Tensor Decompositions} \label{sec:symmetric_tensor_decompositions}

An order-$r$ tensor is a multidimensional array in $\mathbb{K}^{n_1 \times \ldots \times n_r}$, where $\mathbb{K}$ is a field. In this paper, we take $\mathbb{K}=\R$. A tensor $T=(T_{j_1 \ldots j_r}) \in \R^{d_0 \times \ldots \times d_0}$ is symmetric if $T_{j_1 \ldots j_r} = T_{\sigma(j_1) \ldots \sigma(j_r)}$ for all permuations $\sigma \in S_{d_0}$.
Given an order-$r$ symmetric tensor $T=(T_{j_1 \ldots j_r})$ of format $d_0 \times d_0 \times \cdots \times d_0$, one can associate the following homogeneous polynomial of degree $r$ in $d_0$ variables to the tensor $T$:
\begin{equation} \label{eq:homogeneous_polynomial2}
F(\mathbf{x}) = \sum_{1 \leq j_1, j_2, \ldots, j_{r} \leq d_0}  T_{j_1 j_2 \ldots j_r} x_{j_1} x_{j_2} \cdots x_{j_r}. 
\end{equation}
Monomials generally appear more than once in the above sum. 

For a vector $\mathbf{i}=(i_1,\ldots,i_{d_0}) \in \mathbb{N}^{d_0}$ with $i_1+\ldots+i_{d_0}=r$, the multinominal coefficient is defined as
\[
\binom{r}{i_1,\ldots,i_{d_0}} = \frac{r!}{i_1! \cdots i_{d_0}!}.
\] 
Using multinomial coefficients, the polynomial~(\ref{eq:homogeneous_polynomial2}) can be rewritten as
\begin{equation} \label{eq:homogeneous_polynomial}
F(\mathbf{x}) = \sum_{i_1 + i_2 + \ldots + i_{d_0} = r}  \binom{r}{i_1,\ldots,i_{d_0}}   a_{i_1 i_2 \ldots i_{d_0}} x_1^{i_1} x_2^{i_2} \cdots x_{d_0}^{i_{d_0}} 
\end{equation}
such that each monomial appears precisely once in the sum.

To obtain the converse map from homogeneous polynomials to symmetric tensors, we note that there is a bijection between the monomials of degree $r$ in $d_0$ variables and unique entries of a general order-$r$ symmetric tensor of format $d_0 \times d_0 \times \cdots \times d_0$. The bijection is given by the map $f: \mathbb{N}^{d_0} \rightarrow \mathbb{N}^r, \mathbf{i} \mapsto \mathbf{j}$, where $i_k$ denotes the number of appearances of $k$ in $\mathbf{j}$ and the entries of $\mathbf{j}$ are in ascending order. The homogeneous polynomial~(\ref{eq:homogeneous_polynomial})
maps to the symmetric tensor $T$ with the entries
\[
T_{\mathbf{j}} =  a_{f^{-1}(\mathbf{j})},
\]
where the entries of $\mathbf{j}$ are in ascending order. The rest of the entries of $T$ are obtained from the symmetry property.
For more details on the bijection, we refer the reader to~\cite{comon1996decomposition,comon2008symmetric}.

The outer product of $r$ vectors $v_i \in \R^{n_i}$ is an order-$r$ tensor defined as
\[
v_1 \otimes \ldots \otimes v_r = (v_{1,i_1} \cdots v_{r,i_r})_{i_1,\ldots,i_r=1}^{n_1,\ldots,n_r}.
\]

\begin{definition}
Let $T \in \mathbb{K}^{d_0 \times \ldots \times d_0}$ be a symmetric tensor. The \emph{symmetric rank} of $T$ over a field $\mathbb{K}$ is the smallest $k \in \mathbb{N}$ such that
\[
T = \sum_{i=1}^k \lambda_i v_i \otimes \ldots \otimes v_i,
\]
where $\lambda_i \in \mathbb{K}$ and $v_i \in \mathbb{K}^{d_0}$ for $i \in [k]$. If $\mathbb{K}=\R$ (resp. $\mathbb{K}=\C$), then we call it the \emph{real} (resp.\ \emph{complex}) \emph{symmetric tensor rank}.
\end{definition}

Note that for symmetric matrices, the symmetric rank is equal to the rank.

Let $I_1 \cup I_2 = [r]$ be a partition of the set $[r]$. Let $D_j = \prod_{i \in I_j} n_i$ for $j=1,2$. Every partition $I_1 \cup I_2 = [r]$ induces a \emph{flattening} of a tensor $T \in \mathbb{K}^{n_1 \times \ldots \times n_r}$ to a matrix in $\mathbb{K}^{D_1 \times D_2}$. The importance of tensor flattenings is that their minors give relations for tensor rank. A non-zero tensor has rank one if and only if all $2 \times 2$ minors of all its flattenings vanish~\cite[\S 3.4]{landsberg2011tensors}.

\begin{example}
Let $d_0=2$ and $r=3$. Then we consider cubic forms in two variables or, equivalently, order-$3$ $2 \times 2\times 2$ symmetric tensors. The polynomial $x_1^3+3x_1x_2^2+ 3x_2^3$ corresponds to the symmetric tensor whose flattening corresponding to the partition $\{1,2\} \cup \{3\}$ is
$
\begin{pmatrix}
1 & 0 & 0 & 1\\
0 & 1 & 1 & 3 
\end{pmatrix}^T$. 
\end{example}

\subsection{Neuromanifolds for $r=2$} \label{sec:neuromanifolds_r2}

\begin{lemma} \label{lemma:d0d11r2}
The neuromanifold for the architecture $\mathbf{d}=(d_0,d_1,1),r=2$ consists of $(d_0 \times d_0)$ symmetric matrices
\begin{equation} \label{eqn:symmetric_matrix}
\begin{pmatrix}
c_{11} & \frac{1}{2} c_{12} & \cdots & \frac{1}{2} c_{1 d_0} \\
\frac{1}{2} c_{12} & c_{22} & \cdots & \frac{1}{2} c_{2 d_0} \\
\vdots & \vdots & \ddots & \vdots \\
\frac{1}{2} c_{1 d_0} & \frac{1}{2} c_{2 d_0} & \cdots & c_{d_0 d_0} \\
\end{pmatrix}
\end{equation}
of rank at most $d_1$. It is equal to its neurovariety and the ideal of the neurovariety is generated by all $d_1 + 1$ minors of the symmetric matrix. The architecture is filling if and only if $d_1 \geq d_0$. 
\end{lemma}

\begin{proof}
We have $p_{\w}(x) = w_{211}(w_{111} x_1 + \ldots + w_{11 d_0} x_{d_0})^2 + \ldots + w_{21 d_1}(w_{1d_1 1} x_1 + \ldots + w_{1 d_1 d_0} x_{d_0})^2$. Each of the summands of $p_{\w}(x)$ corresponds to a rank-1 symmetric matrix. Since there are $d_1$ summands, the sum corresponds to a symmetric matrix~(\ref{eqn:symmetric_matrix}) of rank at most $d_1$. 
\end{proof}

The next lemma is a special case of~\Cref{lemma:dd01d2_r} and we state it here separately as the special case is used in~\Cref{table:small_examples}.

\begin{lemma} \label{cor:dd01d2_r2}
The neuromanifold for the architecture $\mathbf{d}=(d_0,1,d_2),r=2$ consists of $d_2$ tuples of $d_0 \times d_0$ symmetric matrices of rank at most one such that all the rank-1 matrices are multiples of each other. The neuromanifold is equal to the neurovariety and the ideal of the neurovariety is generated by the $ 2 \times 2$ minors of
\[
\begin{pmatrix}
c_{11}^{(1)} & \frac{1}{2} c_{12}^{(1)} & \cdots & \frac{1}{2} c_{1 d_0}^{(1)} & \cdots & c_{11}^{(d_2)} & \frac{1}{2} c_{12}^{(d_2)} & \cdots & \frac{1}{2} c_{1 d_0}^{(d_2)}  \\
\frac{1}{2} c_{12}^{(1)} & c_{22}^{(1)} & \cdots & \frac{1}{2} c_{2 d_0}^{(1)} & \cdots & \frac{1}{2} c_{12}^{(d_2)} & c_{22}^{(d_2)} & \cdots & \frac{1}{2} c_{2 d_0}^{(d_2)} \\
\vdots & \vdots & \ddots & \vdots & \cdots & \vdots & \vdots & \ddots & \vdots \\ 
\frac{1}{2} c_{1 d_0}^{(1)} & \frac{1}{2} c_{2 d_0}^{(1)} & \cdots & c_{d_0 d_0}^{(1)} & \cdots & \frac{1}{2} c_{1 d_0}^{(d_2)} & \frac{1}{2} c_{2 d_0}^{(d_2)} & \cdots & c_{d_0 d_0}^{(d_2)}
\end{pmatrix}. 
\]
\end{lemma}

The next lemma is due to Maksym Zubkov.

\begin{lemma}
    \label{lem:nn2}
    Let $d_0,d_2 \geq 2$, then for $\mathbf{d}=(d_0,d_0,d_2),r=2$, the parameter map $\Psi_{\mathbf{d},r}$ is not surjective.
\end{lemma}
\begin{proof}
    We will first prove the statement for $\mathbf{d}=(d_0,d_0,2),r=2$. Assume to the contrary that the map $\Psi_{\mathbf{d},r}$ is surjective. For $i \in [d_0]$, we write $\ell_i=(W_1 \x)_i = w_{1i1}x_1+\dots+w_{1id_0}x_{d_0}$. Then for any pair of quadratic forms $(q_1, q_2)\in\Sym_2(\R^{d_0})\times \Sym_2(\R^{d_0})$, we have some $\w$ such that
    \[\scalebox{0.9}{$\Psi_{\mathbf{d},r}(\w)=\begin{pmatrix}
        w_{211}\ell_1^2+\dots+w_{21 d_0}\ell_{d_0}^2\\
        w_{221}\ell_1^2+\dots+w_{22 d_0}\ell_{d_0}^2\\
    \end{pmatrix}=\begin{pmatrix}
        \ell_1 w_{211} \ell_1+\dots+\ell_{d_0} w_{21d_0} \ell_{d_0}\\
        \ell_1 w_{221} \ell_1+\dots+\ell_{d_0} w_{22d_0} \ell_{d_0}\\
    \end{pmatrix}
    =\begin{pmatrix}
        \x^T W_1^T D_1 W_1 \x\\
        \x^T W_1^T D_2 W_1 \x
    \end{pmatrix}$}\]
    where $D_i$ is the diagonal matrix with $w_{2ij}$ for $i=1,2$ and $j=1,\dots,d_0$. So, if $A_1$ and $A_2$ are the corresponding matrices for quadratic forms $q_1$ and $q_2$, then we would have
    \[A_1=W_1^T D_1 W_1,\quad A_2=W_1^T D_2 W_1.\]
    This tells us that the two symmetric matrices $A_1$ and $A_2$ can be diagonalized simultaneously which is in general not true as $A_1$ and $A_2$ do not commute in general for $d_0\geq 2$. Therefore, $\Psi_{\d,r}$ is not surjective. 

    The proof can be extended to the architectures $(d_0,d_0,d_2)$ for any $d_2\geq 2$ as it corresponds to simultaneous diagonalization of $d_2$ matrices which is impossible for an arbitrary choice of $d_2$ matrices. 
\end{proof}

\subsection{Neuromanifolds for $r\in\N$} \label{sec:neuromanifolds_r>2}

\begin{lemma} \label{lemma:dd01d2_r}
The neuromanifold for the architecture $\mathbf{d}=(d_0,1,d_2),r \in \N$ consists of $d_2$ order-$r$ $(d_0 \times d_0 \times \cdots \times d_0)$ symmetric tensors of rank at most one such that all the rank-1 tensors are multiples of each other. The neuromanifold is equal to the neurovariety and the ideal of the neurovariety is generated by the $2 \times 2$ minors of the flattenings of the $(d_0  \times \cdots \times d_0 \times d_0 d_2)$ tensor that is obtained from combining the $d_2$ $(d_0 \times d_0 \times \cdots \times d_0)$ symmetric tensors along the last index.
\end{lemma}

\begin{proof}
We have 
\[
p_{\w}(x) = 
\begin{pmatrix}
w_{211}(w_{111} x_1 + \ldots + w_{11 d_0} x_{d_0})^r\\
w_{221}(w_{111} x_1 + \ldots + w_{11 d_0} x_{d_0})^r\\
\vdots \\
w_{2d_2 1}(w_{111} x_1 + \ldots + w_{11 d_0} x_{d_0})^r
\end{pmatrix}.
\]
The form $(w_{111} x_1 + \ldots + w_{11 d_0} x_{d_0})^r$ corresponds to an order-$r$ rank-1 symmetric tensor. The $i^{\text{th}}$ component of $p_{\w}(x)$ is the symmetric tensor multiplied by the constant $w_{2 i 1}$. This proves the first statement. For the second statement, we note that combining the rank-1 tensors to a $d_0  \times \cdots \times d_0 \times d_0 d_2$ tensor along the last index gives another rank-1 tensor. Moreover, any $d_0  \times \cdots \times d_0 \times d_0 d_2$ tensor of rank one such that each of the $d_2$  $(d_0 \times d_0 \times \cdots \times d_0)$-subtensors is symmetric and can be obtained from symmetric $(d_0 \times d_0 \times \cdots \times d_0)$ rank-1 tensors that are multiples of each other. The second statement follows from the result that the ideal of the tensors of rank at most one is generated by $2\times 2$ minors of flattenings~\cite[\S 3.4]{landsberg2011tensors}.
\end{proof}

\begin{proposition}
    \label{prop:fillingMfd}
    Let $\d = (d_0, d_1, d_2)$ with $d_1 \geq \binom{r+d_0-1}{r}$. Then the neuromanifold itself is filling, i.e.,\ $\M_{\d,r} = (\Sym_r(\R^{d_0}))^{d_2}$.
\end{proposition}

\begin{proof}
    Let $N\coloneqq \binom{r+d_0-1}{r}$.
    Write the coefficients of $p_{\w} = \Psi_{\d,r}(\w)$ as
    \[
        C = \begin{pmatrix}
        (c_{I}^{(1)})_I\\
        \vdots\\
        (c_{I}^{(d_2)})_I
        \end{pmatrix}=\begin{pmatrix}
            w_{211} & w_{212} & \hdots & w_{21d_1}\\ \vdots & \vdots & \ddots & \vdots \\ 
            w_{2d_2 1} & w_{2d_2 2} & \hdots & w_{2d_2d_1}
        \end{pmatrix}
        \begin{pmatrix}
            (W_{1,1})^I \\
            \vdots \\
            (W_{1,d_1})^I
        \end{pmatrix},
    \]
    where $I$ ranges over all multiindices in $\binom{[d_0]}{r}$ and $(W_{1,i})^I$ denotes the row consisting of all degree-$r$ monomials formed by entries of the $i^{\text{th}}$ row of $W_1$ with the corresponding multinomial coefficients.
    We want to show that any real $d_2\times N$ matrix $C$ can be written in this form.
    If $d_1 \geq N$ then for a general choice of parameters $W_1$ this matrix has a left-inverse with real entries and hence the system above can be solved over the reals.
\end{proof}

\begin{definition}
    Let $\mathbb{k}$ be $\R$ or $\C$ and let 
    \[
        S^n_{d,k}(\mathbb{k}) = \{T \in \Sym_d(\mathbb{k}^n): \text{ symmetric rank} \\ \text{ of } T \text{ is } k\}.
    \] 
    If $S^n_{d,k}(\mathbb{k}) \subseteq \Sym_d(\mathbb{k}^n)$ has non-empty interior with respect to the Euclidean topology, $k$ is called a \emph{typical symmetric rank} for order-$d$ tensors of dimension $n$.
\end{definition}

Over the complex numbers, there is a unique typical symmetric rank which is called the generic symmetric rank. However, over the real numbers there can be multiple typical symmetric ranks. In the following we will write $S^n_{d,k}$ for $S^n_{d,k}(\R)$. A neuromanifold $\M_{\d,r}$ having non-empty interior with respect to the Euclidean topology implies that the architecture $(\d,r)$ is filling.  The importance of typical symmetric ranks for PNNs comes from the following fact which is a consequence of the discussions above and the inclusion $S^n_{d,k} \subseteq \M_{(n,k,1),d}$.

\begin{theorem}
    \label{thm:typicalRanksNeuro}
    If $k$ is a typical symmetric rank for order-$d$ tensors of dimension $n$ then the architecture $\d=(n,k,1),~r=d$ is filling, i.e., $\V_{\d,r} = \Sym_d(\R^n)$. Moreover, if there are two typical symmetric ranks $k < k'$ then $\mathrm{cl_{Eucl}}(\M_{(n,k,1),d}) \subsetneq \Sym_d(\R^n)$ where $\mathrm{cl_{Eucl}}$ denotes the closure with respect to the Euclidean topology.
\end{theorem}

Let us consider the case $n=2$ first; this is the setting of the paper \cite{comon2012typical}. We summarize their findings below.

\begin{theorem}[{\cite[Proposition 2.2 \& Main Theorem]{comon2012typical}}]
    \hfill
    \begin{enumerate}
        \item $S^2_{3,k}$ has non-empty interior only for $k \in \{2,3\}$.
        \item $S^2_{4,k}$ has non-empty interior only for $k \in \{3,4\}$.
        \item $S^2_{5,k}$ has non-empty interior only for $k \in \{3,4,5\}$.
    \end{enumerate}
\end{theorem}

Using this result we obtain the following.

\begin{corollary} \label{cor:neuromanifolds_from_tensors1}
    \hfill
    \begin{enumerate}
        \item $\d = (2,d_1,1),~r=3$ is filling for $d_1\geq 2$ and 
        \[\mathrm{cl_{Eucl}}(\M_{(2,2,1),3}) \subsetneq \mathrm{cl_{Eucl}}(\M_{(2,3,1),3}) = \Sym_{3}(\R^2).\]
        \item $\d = (2,d_1,1),~r=4$ is filling for $d_1\geq 3$ and 
        \[\mathrm{cl_{Eucl}}(\M_{(2,3,1),4}) \subsetneq \mathrm{cl_{Eucl}}(\M_{(2,4,1),4}) = \Sym_{4}(\R^2).\]
        \item $\d = (2,d_1,1),~r=5$ is filling for $d_1\geq 3$ and 
        \[\mathrm{cl_{Eucl}}(\M_{(2,3,1),5}) \subsetneq \mathrm{cl_{Eucl}}(\M_{(2,4,1),5}) \subsetneq \mathrm{cl_{Eucl}}(\M_{(2,5,1),5}) = \Sym_{3}(\R^2).\]
    \end{enumerate}
\end{corollary}

Some similar results on typical symmetric ranks are also known for higher dimensional tensors though the difficulty increases quickly. We summarize some consequences below.

\begin{corollary} \label{cor:neuromanifolds_from_tensors2}
    \hfill
    \begin{enumerate}
        \item $\d = (3,d_1,1),~r=4$ is filling for $d_1\geq 6$ and 
        \[\mathrm{cl_{Eucl}}(\M_{(3,6,1),4}) \subsetneq \mathrm{cl_{Eucl}}(\M_{(3,7,1),4}) \subseteq \mathrm{cl_{Eucl}}(\M_{(3,8,1),4}) = \Sym_{4}(\R^3).\]
        \item $\d = (3,d_1,1),~r=5$ is filling for $d_1\geq 7$ and 
        \[\mathrm{cl_{Eucl}}(\M_{(3,7,1),5}) \subsetneq \mathrm{cl_{Eucl}}(\M_{(3,8,1),5}) \subseteq \dots \subseteq \mathrm{cl_{Eucl}}(\M_{(3,13,1),5}) = \Sym_{5}(\R^3).\]
        \item $\d = (4,d_1,1),~r=3$ is filling for $d_1\geq 5$ and 
        \[\mathrm{cl_{Eucl}}(\M_{(4,5,1),3}) \subsetneq \mathrm{cl_{Eucl}}(\M_{(4,6,1),3}) = \Sym_{3}(\R^4)\]
    \end{enumerate}
\end{corollary}

\begin{proof}
    The statements follow from Theorem \ref{thm:typicalRanksNeuro} and \cite[Theorem 1.2]{realTypicalRanks}.
\end{proof}

The algebraic boundary $\del_{\text{alg}}(\M)$ of a set $\M$ is the Zariski closure of its topological boundary $\del(\M)$. Little is known about the algebraic boundaries in the cases considered above. Theorem 4.1 in \cite{Michaek2016} implies that the algebraic boundary $\del_{\text{alg}}(\M_{(3,6,1),4})$ has an irreducible component of degree 51.

\section{Neurovarieties} \label{sec:neurovarities}

In this section, we will describe different approaches to studying neurovarieties. We start by studying the fibers of a neurovariety to characterize the neurovariety $\V_{(2,2,d_2),2}$ (\Cref{prop:d22k_r2}) and partially the neuromanifold $\M_{(2,2,d_2),2}$ (\Cref{prop:semialgebraic_description_d22k_r2}). We use Grassmannians to describe the neurovariety $\V_{(3,3,3),2}$ (\Cref{example:grassm}) and apply the Hilbert--Burch Theorem to study $\V_{(2,2,2,2),2}$ (\Cref{example:Hilbert-Burch}).

\begin{proposition} \label{prop:d22k_r2}
    The neurovariety $\V_{(2,2,d_2),2}$ is the vanishing locus of all $3\times 3$ minors of \[C_{d_2} = \begin{pmatrix}
    c_{11}^{(1)} & c_{12}^{(1)} & c_{22}^{(1)} \\
    \vdots & \vdots & \vdots \\
    c_{11}^{(d_2)} & c_{12}^{(d_2)} & c_{22}^{(d_2)}
\end{pmatrix}=\begin{pmatrix}
    w_{211} & w_{212} \\ \vdots & \vdots \\ w_{2d_2 1} & w_{2d_2 2}
\end{pmatrix}
\begin{pmatrix}
    w_{111}^2 & 2w_{111}w_{112} & w_{112}^2 \\
    w_{121}^2 & 2w_{121}w_{122} & w_{122}^2
\end{pmatrix}.\]  
\end{proposition}
\begin{proof}
Clearly $\V_{(2,2,d_2),2} \subseteq \V(\langle 3\times 3\text{ minors of } C_{d_2}\rangle)$ as $\rk(C_{d_2}) \leq 2$. Now, $\V(\langle 3\times 3\text{ minors of } C_{d_2}\rangle)$ is irreducible and $\dim(\V(\langle 3\times 3\text{ minors of } C_{d_2}\rangle)) = 2d_2 + 2$, see e.g.\ \cite[Proposition 1.1]{bruns2006determinantal}. Consider the parameter map and the corresponding neurovariety
\[\Psi_{(2,2,d_2),2} \colon \underbrace{\R^{2\times 2}\times \R^{d_2\times 2}}_{\dim 2d_2+4} \rightarrow (\Sym_{2}(\R^2))^{d_2}, \quad \V_{(2,2,d_2),2} = \overline{\im(\Psi_{(2,2,d_2),2})}.\] 
The fibers of $\Psi_{(2,2,d_2),2}$ over $\V_{(2,2,d_2),2}$ are at least $2$-dimensional by Lemma \ref{lem:multiHomog}. 
For $d_2=2$, one can computationally check that over a generic matrix $C_{d_2}$ the fiber is 2-dimensional. Since generic fibers are at least 2-dimensional, they are exactly 2-dimensional. 
For the case $d_2 > 2$, note that by setting the parameters $w_{211},w_{212},\dots,w_{2(d_2-1)1}, w_{2(d_2-1)2}$ to zero and varying $w_{2d_2 1}, w_{2d_2 2}$, we obtain a 2-dimensional family not contained in the embedding of $\M_{(2,2,d_2-1),2}$ into $\M_{(2,2,d_2),2}$. As the dimension of the parameter space increases by two when passing from $d_2-1$ to $d_2$, the generic fiber is 2-dimensional by induction.
Thus $\dim(\V_{(2,2,d_2),2}) = 2d_2+2$ and moreover $\V_{(2,2,d_2),2}$ is irreducible, hence $\V_{(2,2,d_2),2} = \V(\langle 3\times 3\text{ minors of } C_{d_2}\rangle))$.
\end{proof}

\begin{proposition} \label{prop:semialgebraic_description_d22k_r2}
    For $\d = (2,2,d_2)$ ($d_2\geq 2$) and $r=2$, the Euclidean closure of $\M_{\d,r}$ is strictly included in $\V_{\d,r}$. In the case $\d=(2,2,2)$, we have 
    \[
        C = \begin{pmatrix}
            c_{11}^{(1)} & c_{12}^{(1)} & c_{22}^{(1)} \\
            c_{11}^{(2)} & c_{12}^{(2)} & c_{22}^{(2)}
            \end{pmatrix}
        \in \M_{(2,2,2),2}
    \]
    if and only if 
    \[
        M_{1,3}(C)^2 \geq M_{1,2}(C)M_{2,3}(C),
    \]
    where $M_{i,j}(C)$ is the $2\times 2$ minor of $C$ obtained by taking the determinant of the $i^{\text{th}}$ and $j^{\text{th}}$ column.
    In particular, the algebraic boundary of $\M_{(2,2,2),2}$ is given by
    \[
        \del_{\text{alg}}\M_{(2,2,2),2} = 
        \V(M_{1,3}(C)^2 - M_{1,2}(C)M_{2,3}(C)).
    \]
\end{proposition}

\begin{proof}
    First consider the case $\d = (2,2,2)$. We are looking for conditions on a matrix $A\in \R^{2\times 3}$ such that there exists $G\in GL_2(\R)$ with 
    \begin{equation}
        \label{equ:GAW}
        \begin{pmatrix}
            g_{11} & g_{12} \\
            g_{21} & g_{22} \\
        \end{pmatrix}
        \begin{pmatrix}
             a_{11} & a_{12} &  a_{13} \\
             a_{21} & a_{22} &  a_{23} 
        \end{pmatrix}
        = \begin{pmatrix}
        w_{111}^2 & 2w_{111}w_{112} & w_{112}^2 \\
        w_{121}^2 & 2w_{121}w_{122} & w_{122}^2
    \end{pmatrix}.
    \end{equation}
     We consider two cases: when both $w_{111},w_{121}$ are non-zero and when at least one of them is zero. If $w_{111}$ and $w_{121}$ are non-zero, then after rescaling, w.l.o.g.\ we can assume $w_{111}=w_{121}=1$. Each matrix entry in \eqref{equ:GAW} yields a polynomial equation; these generate an ideal $I$ in the polynomial ring 
     $\R[g_{11},\dots, g_{22}, w_{112}, w_{122}, a_{11}, \dots, a_{23}].$ 
     Two generators in a Gr\"obner basis of $I$ with respect to lexicographic ordering are
    \begin{equation} \label{eq:222_GB_generator}
        w_{1i2}^2a_{11}a_{22}-w_{1i2}^2a_{12}a_{21}
        -2w_{1i2}a_{11}a_{23}+2w_{1i2}a_{13}a_{21}+a_{12}a_{23}-a_{13}a_{22},
    \end{equation}
    where $i\in\{1,2\}$. These quadratic equations in $w_{1i2}$ have discriminant
    \[
        \Delta = M_{1,3}(A)^2 - M_{1,2}(A)M_{2,3}(A),
    \]
    where $A$ is the $2\times 3$ matrix with entries $a_{ij}$ from \eqref{equ:GAW}. The denominator of the quadratic formula for~\eqref{eq:222_GB_generator} is $M_{1,2}(A)$. Therefore, real solutions to \eqref{equ:GAW} with $w_{111},w_{121}$ non-zero exist if and only if $M_{1,3}(A)^2 \geq M_{1,2}(A)M_{2,3}(A)$ and $M_{1,2}(A) \neq 0$. \par 

    When at least one $w_{1i1}=0$, then possibly after rescaling, the Gr\"obner basis of $I$ with respect to the lexicographic ordering contains $w_{i12}^2 a_{11} a_{22} - w_{i12}^2 a_{12} a_{21}$. This implies $w_{i12}=0$ or $M_{1,2}(A)=0$.  The former condition implies $M_{1,2}(A)=0$, so we do not have to consider it separately. When $M_{1,2}(A)=0$, the inequality $M_{1,3}(A)^2 \geq M_{1,2}(A)M_{2,3}(A)$ is automatically satisfied. Therefore, real solutions to \eqref{equ:GAW} with at least one of $w_{111},w_{121}$ zero exist if and only if $M_{1,3}(A)^2 \geq M_{1,2}(A)M_{2,3}(A)$ and $M_{1,2}(A) = 0$. Combining the two cases gives that real solutions to \eqref{equ:GAW} exist if and only if $M_{1,3}(A)^2 \geq M_{1,2}(A)M_{2,3}(A)$.

    As this inequality remains unchanged after multiplying $A$ with a $2\times 2$ matrix (both sides get multiplied by the squared determinant of the matrix which is positive), we conclude $C\in \M_{(2,2,2),2}$ if and only if 
    $M_{1,3}(C)^2 \geq M_{1,2}(C)M_{2,3}(C)$. 

    In the case $d_2>2$, it is clear that at least for any two rows $i,j$ of $C$, the inequality 
    \[
        M_{1,3}(C_{i,j})^2 \geq M_{1,2}(C_{i,j})M_{2,3}(C_{i,j})
    \]
    needs to hold. Thus we get a full-dimensional set of points in $\V_{\d,r}\setminus \M_{\d,r}$, hence $\mathrm{cl_{Eucl}}(\M_{\d,r}) \subsetneq \V_{\d,r}$.
\end{proof}

\begin{remark}
    \label{rem:(d_0,2,d_2),r=2}
    For example, if $A =     \begin{pmatrix}
         a & \star &  -a \\
         b & \star &  -b 
    \end{pmatrix}$ in the proof of~\Cref{prop:semialgebraic_description_d22k_r2}, then the equation
\[
    \begin{pmatrix}
        g_{11} & g_{12} \\
        g_{21} & g_{22} \\
    \end{pmatrix}
    \begin{pmatrix}
         a & \star &  -a \\
         b & \star &  -b 
    \end{pmatrix}
    = \begin{pmatrix}
        w_{111}^2 & 2w_{111}w_{112} & w_{112}^2 \\
        w_{121}^2 & 2w_{121}w_{122} & w_{122}^2
    \end{pmatrix}
    \]
does not have a solution.
    This example can be extended to any architecture $\d=(d_0,2,d_2)$ with $d_0,d_2\geq 2$ and activation degree $r=2$ to show that the neuromanifold $\M_{\d,r}$ is not equal to $\V_{\d,r}$.
\end{remark}

In the rest of the section, we use more advanced tools from algebraic geometry to compute neurovarieties.\ A na\"ive elimination approach fails to compute these varieties.

\begin{example} 
\label{example:grassm}
Consider the architecture $\d=(3,3,3),~r=2$; the image of $\Psi_{\mathbf{d},r}$ consists of three linear combinations of three squares of linear forms in three variables.
Let us take a $3\times 6$ matrix representing the embedding of $\P^2_x\times\P^2_y\times\P^2_z$ via the degree two Veronese map (with corresponding scaling) into $\P^5\times\P^5\times\P^5$:
\[
    \begin{pmatrix}
        x_0^2 & 2x_0x_1 & 2x_0x_2 & x_1^2 & 2x_1x_2 & x_2^2 \\
        y_0^2 & 2y_0y_1 & 2y_0y_2 & y_1^2 & 2y_1y_2 & y_2^2 \\
        z_0^2 & 2z_0z_1 & 2z_0z_2 & z_1^2 & 2z_1z_2 & z_2^2 \\
    \end{pmatrix}    
\]
By taking all the $3\times 3$ minors (20 in total) of this matrix, we get an embedding into the Grassmannian $\gr(3,6)$ in its Pl\"ucker coordinates.
The ideal of this variety in Pl\"ucker coordinates can be found using elimination. Its dimension is six and its degree is 57.
Now for each of the 20 Pl\"ucker coordinates we substitute a $3\times 3$ minor of a $3\times 6$ matrix of unknowns, where each unknown represents a coefficient of a degree two monomial (we view $\P^5$ as the space of all homogeneous quadrics in three variables):
\[
    C = \begin{pmatrix}
        c_{11}^{(1)} & c_{12}^{(1)} & \dots & c_{33}^{(1)} \\
        c_{11}^{(2)} & c_{12}^{(2)} & \dots & c_{33}^{(2)} \\
        c_{11}^{(3)} & c_{12}^{(3)} & \dots & c_{33}^{(3)} \\
    \end{pmatrix}
\]
The resulting ideal $I$ gives rise to a variety $X=\V(I)$. This variety has two irreducible components: one is given by the vanishing of all $3\times 3$ minors of $C$, the other is the neurovariety $\V_{(3,3,3),2}$ embedded into $\P^5_{c^{(1)}}\times\P^5_{c^{(2)}}\times\P^5_{c^{(3)}}$.
The codimension of $\V_{(3,3,3),2}\subseteq \P^5_{c^{(1)}}\times\P^5_{c^{(2)}}\times\P^5_{c^{(3)}}$ is three, hence the affine dimension is $\dim(\V_{\d,r}) = 18-3=15$. \par
The first part of the computation can be found in \cite{brodsky2011tropical}; we implemented the whole procedure and made it available at \cite{mathrepo}.
\end{example}

We would like to thank Bernd Sturmfels for telling us about the following example; the computation is motivated by the Hilbert--Burch Theorem, see for instance \cite[Theorem 20.15]{eisenbud2013commutative}.

\begin{theorem}[Hilbert--Burch]
    Let $R$ be a local ring with an ideal $I\subset R$ such that 
    \[
        R^n \xrightarrow{f} R^m \rightarrow R \rightarrow R/I \rightarrow 0
    \]
    is a free resolution of $R/I$, then $m=n-1$ and $I=aJ$ where $a$ is a regular element of $R$ and $J$ is a depth-two ideal generated by all $m\times m$ minors of the matrix representing $f$.
\end{theorem}

Consider the architecture $(\d = (2,2,2,2),r=2)$. We expect that $\V_{\d,r}$ has dimension eight and thus has codimension two in $(\Sym_{4}(\R^2))^2$. Therefore, we can hope that by the Hilbert--Burch Theorem we can find a matrix whose minors generate the ideal of $\V_{\d,r}$. Again, we work in projective space, i.e., we consider $\V_{\d,r}$ as a projective variety inside $\P^4_{c^{(1)}} \times \P^4_{c^{(2)}}$.

\begin{example} \label{example:Hilbert-Burch}
One finds that the ideal cutting out $\V_{(2,2,2,2),2}\subseteq \P^4_{c^{(1)}} \times \P^4_{c^{(2)}}$ is generated by the $5\times 5$ minors of the following matrix
\[
    \begin{pmatrix*}[r]
        c_{1112}^{(2)} & -c_{1111}^{(2)} & \makebox[\widthof{$c_{1112}^{(2)}$}][c]{0} & c_{1112}^{(1)} & -c_{1111}^{(1)} & \makebox[\widthof{$c_{1112}^{(2)}$}][c]{0} \\
        4c_{1122}^{(2)} & -c_{1112}^{(2)} & 8c_{1111}^{(2)} & 4c_{1122}^{(1)} & -c_{1112}^{(1)} & 8c_{1111}^{(1)} \\
        8c_{1222}^{(2)} & \makebox[\widthof{$c_{1112}^{(2)}$}][c]{0} & 8c_{1112}^{(2)} & 8c_{1222}^{(1)} & \makebox[\widthof{$c_{1112}^{(2)}$}][c]{0} & 8c_{1112}^{(1)} \\
        8c_{2222}^{(2)} & c_{1222}^{(2)} & 4c_{1122}^{(2)} & 8c_{2222}^{(1)} & c_{1222}^{(1)} & 4c_{1122}^{(1)} \\
        \makebox[\widthof{$c_{1112}^{(2)}$}][c]{0} & c_{2222}^{(2)} & c_{1222}^{(2)} & \makebox[\widthof{$c_{1112}^{(2)}$}][c]{0} & c_{2222}^{(1)} & c_{1222}^{(1)} \\
    \end{pmatrix*}.  
\]
\end{example}

\section{Dimension} \label{sec:dimension}

\subsection{A plethora of conjectures}
Determining dimensions of neurovarieties is a difficult task and there exists a wide range of open conjectures in this direction. A classical result is the Alexander--Hirschowitz Theorem; already conjectured in the late $19^{\text{th}}$ century, the proof was completed about a century later in 1995. Here is a formulation in the language of PNNs.

\begin{theorem}[Alexander--Hirschowitz \cite{AH}]
    \label{thm:AH}
    If $\d = (d_0,d_1,1)$, $\V_{\d,r}$ attains the expected dimension $\min\{d_0d_1, \binom{d_0+r-1}{r}\}$, except for the following cases:
    \begin{enumerate}
        \item $r=2,~2\leq d_1 < d_0$,
        \item $r=3,~d_0=5,~d_1=7$ where $\dim(\V_{(5,7,1),3}) = 34$ (defect 1),
        \item $r=4,~d_0=3,~d_1=5$ where $\dim(\V_{(3,5,1),4}) = 14$ (defect 1),
        \item $r=4,~d_0=4,~d_1=9$ where $\dim(\V_{(4,9,1),4}) = 34$ (defect 1),
        \item $r=4,~d_0=5,~d_1=14$ where $\dim(\V_{(5,14,1),4}) = 69$ (defect 1).
    \end{enumerate}
\end{theorem}

It is conjectured that the neurovariety attains the expected dimension for large activation degree. Compare this for example with the Alexander--Hirschowitz Theorem which tells us there exist no defective single-output two-layer neurovarieties with activation degree $r\geq 5$.   The  conjecture appeared as a theorem in~\cite{kileel2019expressive}, but a mistake in the proof was pointed out by Theo Elenius in his BSc thesis. After the present article first appeared as a preprint on arXiv, the conjecture has been proven in~\cite[Theorem 12]{finkel2024activation}.

\begin{conjecture}[\cite{kileel2019expressive}, Theorem 14]
    \label{conj:asympEDim}
    For any fixed widths $\mathbf{d}=(d_0,\dots,d_L)$ with $d_i>1$ for $i=1,\dots,L-1$, there exists $\Tilde{r} = \Tilde{r}(\d)$ such that whenever $r>\Tilde{r}$, the neurovariety $\V_{\mathbf{d},r}$ attains the expected dimension.
\end{conjecture}

It is, however, not true that $\V_{\d,r}$ being defective implies that $\V_{\d,r'}$ is defective for all $r'<r$. For example, $\V_{(5,7,1),3}$ has defect one whereas $\V_{(5,7,1),2}$ is non-defective.\par 

Also note that the assumption $d_i>1$ for $i=1,\dots,L-1$ is necessary as is shown by the following example.

\begin{example}
    Consider the widths $\mathbf{d}=(2,1,2,1)$. We claim that for any $r>1$, $\V_{\mathbf{d},r}$ has defect one, i.e., $\dim(\V_{\mathbf{d},r}) = 2$. Indeed, observe that the parameter map $\Psi_{\mathbf{d},r}$ is given by
    \begin{align*}
        \x & \mapsto (w_{111}x_1 + w_{112}x_2)^r \mapsto \begin{pmatrix}
            w_{211}^r(w_{111}x_1 + w_{112}x_2)^{r^2} \\
            w_{221}^r(w_{111}x_1 + w_{112}x_2)^{r^2} \\
        \end{pmatrix} \\
        & \mapsto (w_{311}w_{211}^r + w_{312}w_{221}^r)(w_{111}x_1 + w_{112}x_2)^{r^2}.
    \end{align*}
    As all weights coming from the last two layers can be factored out, we get that $\V_{\mathbf{d},r} \cong \V_{(2,1,1,1),r}$. The latter neurovariety is immediately seen to have dimension two. This generalizes to the following statement.
\end{example}

\begin{proposition} \label{prop:neuromanifold-with-width-1-layer}
    Let $\mathbf{d}_0 \in \mathbb{N}^{n_0}$ and $\mathbf{d}_2 \in \mathbb{N}^{n_2}$. For $\mathbf{d}=(\mathbf{d}_0,1,\mathbf{d}_2,d_3)$, the neuromanifold $\M_{\mathbf{d},r}$ is equal to $\M_{(\mathbf{d}_0,1,\mathbf{1},d_3),r}$ for any $r$. In particular, for $d_3=1$, $\dim(\M_{\mathbf{d},r}) = \dim(\M_{(\mathbf{d}_0,1),r})$.
\end{proposition}

\begin{proof}
Fix parameters $\w=(W_1,W_2,\dots,W_L)$. Let $p = \Psi_{(\d_0,1),r}(\pi_1(\w))$ and $q = \Psi_{(1,\d_2,1)}(\pi_2(\w))$, where $\pi_1$ denotes the projection to the parameter space corresponding to the first $n_0+1$ layers, and $\pi_2$ denotes the projection to the remaining parameters. Then $p_{\w}(\x)$ factors as a composition $p_{\w}(\x) = (q\circ p)(\x)$. The polynomial $q$ is a single-input neural network, hence we can write $(q\circ p)_i(\x)$ as a product $\alpha_i(\pi_2(\w))\cdot p(\x)^{r^{n_2}}$ where $\alpha_i(\pi_2(\w))$ only depends on weights corresponding to the last $n_2+2$ layers and thus 
\[
    \M_{\mathbf{d},r} = \{(\alpha_i \cdot p^{r^{n_2}})_i \colon \alpha_i \in \mathbb{R}, p \in \M_{(\mathbf{d}_0,1),r}\}.
\]
The latter set is by definition equal to  $\M_{(\mathbf{d}_0,1,\mathbf{1},d_3),r}$. In the case $d_3=1$, its dimension is equal to $\dim(\M_{(\mathbf{d}_0,1),r})$.
\end{proof}

\begin{remark}
This allows to give a strong bound on the dimension of $\M_{\d,r}$ for $\d=(\d_0,1,\d_2,d_3)$:
\[
    \dim(\M_{\d,r}) = \dim(\M_{\d_0,1,\mathbf{1},d_3},r) \leq \mathrm{edim}(\M_{(\d_0,1,\mathbf{1},d_3),r}) = \dim(\M_{\d_0,r})+d_3-1.
\]
If $d_0>1$, this is strictly less than the ambient dimension which is $d_3\binom{r^{L-1} + d_0 - 1}{r^{L-1}}$, where $d_0$ is the width of the first layer. Hence the width one at the $(n_0+1)^{\text{st}}$ layer is an instance of an \emph{asymptotic bottleneck}, i.e., the architecture is not filling regardless of adding more layers or changing the widths in all other layers except the input and $(n_0+1)^{\text{st}}$ layer \cite[Def.\ 18]{kileel2019expressive}.
\end{remark}

\begin{remark}
 For $\mathbf{d}=(d_0,1,d_2,1)$, the neuromanifold $\M_{\mathbf{d},r}$ is equal to $\M_{(d_0,1,1,1),r}$ for any $r$. In particular, $\dim(\M_{\mathbf{d},r}) = d_0$.
    This enables us to compute the dimension of many defective neurovarieties. For example, $\dim(\V_{(3,1,5,1),3}) = 3$ and therefore $\V_{(3,1,5,1),3}$ has defect 4.
\end{remark}

In contrast to the asymptotic statement for large activation degree in Conjecture \ref{conj:asympEDim}, we also conjecture the following.

\begin{conjecture}
    Let $\d=(d_0,d_1,\dots,d_L)$ be a non-increasing sequence with $d_L > 1$. Then for any $r$, the neurovariety $\V_{\d,r}$ attains the expected dimension.
\end{conjecture}

We verified this conjecture using Algorithm \ref{alg:backpropagation} described in \S\ref{sec:backprop} for four- and five-layer PNNs with widths up to three and activation degree up to five, see \cite{mathrepo}. In the context of feedforward neural networks with ReLU activation this statement has been proved \cite[Theorem 8.11]{grigsby2022functional}. \par 
Using Lemma \ref{lem:multiHomog}, one can deduce the following recursive dimension bound.

\begin{proposition}[Recursive bound, {\cite[Proposition 17]{kileel2019expressive}}]\label{prop:recursive}
    For any $r\in\N$, $\d = ( d_0,\dots,d_i,\dots,d_L )$ and $i \in \{1,\ldots,L-1\}$, we have
    \[
    \dim (\V_{\d,r}) \leq \dim (\V_{(d_0,\ldots,d_i),r}) + \dim (\V_{(d_i,\ldots,d_L),r}) - d_i. 
    \]
\end{proposition}

\begin{corollary} \label{cor:defective-subnetwork-implies-defective-full-network}
    Let $\tilde{\d} = (d_0,\dots,d_{i-1},\d,d_{j+1},\dots,d_L)$, where $\d=(d_i,d_{i+1},\dots,d_j)$ is an architecture such that $\V_{\d,r}$ is defective for some $r>0$. Here we allow $i=0$ or $j=L$. Suppose $\mathcal{V}_{\tilde{\d},r}$ is not filling, then $\mathcal{V}_{\tilde{\d},r}$ is defective.
\end{corollary}

\begin{proof}
It follows from Proposition \ref{prop:recursive} that 
\begin{equation}\label{eq:dim-add}
    \dim (\V_{\tilde{\d},r}) \leq \dim (\V_{(d_0,\ldots,d_i),r}) +\dim (\V_{\d,r})+ \dim (\V_{(d_j,\ldots,d_L),r}) - d_i-d_j
\end{equation}
where $\dim (\V_{(d_0,\ldots,d_i),r})$ or $\dim (\V_{(d_j,\ldots,d_L),r})$ might not appear in the sum above.
By \Cref{def:expected-dim}, the expected dimension is defined as 
\[\edim(\V_{\d,r}) = \min\left\{d_L + \sum_{i=0}^{L-1}(d_id_{i+1} - d_{i+1}), d_L\binom{d_0 + r^{L-1} - 1}{r^{L-1}}\right\}.\]
Since $\dim (\V_{\d,r})<\edim(\V_{\d,r})$, it follows from equation \eqref{eq:dim-add} that 
\begin{align*}
    \dim (\V_{\tilde{\d},r}) &< d_L + \sum_{\alpha=0}^{i-1}(d_\alpha d_{\alpha+1} - d_{\alpha+1}) + \sum_{\beta=i}^{j-1}(d_\beta d_{\beta+1} - d_{\beta+1}) + \sum_{\gamma=j}^{L}(d_\gamma d_{\gamma+1} - d_{\gamma+1}) \\
    &= d_L + \sum_{\alpha=0}^{L-1}(d_\alpha d_{\alpha+1} - d_{\alpha+1})
\end{align*}
where $\sum_{\alpha=0}^{i-1}(d_\alpha d_{\alpha+1} - d_{\alpha+1})$ or $\sum_{\gamma=j}^{L}(d_\gamma d_{\gamma+1} - d_{\gamma+1})$ might not appear in the sum above.
Since $\V_{\tilde{\d},r}$ is non-filling, it follows directly that $\V_{\tilde{\d},r}$ is defective. 
\end{proof}

\begin{remark}
    The converse to this statement is not true. For example, consider $\d = (2,2,1,2),~r=2$. Then $\V_{\d,r}$ has defect one, but both $\V_{(2,2,1),2}$ and $\V_{(2,1,2),2}$ are non-defective.
\end{remark}

Let $\preccurlyeq$ be the partial order on the set of widths induced by coordinatewise comparison. The following conjecture suggests a unimodal distribution of layer widths within a neural network is efficient. In the realm of machine learning theory, this architectural pattern is commonly believed to allow the network to initially expand its capacity for feature representation, enabling the network to capture intricate patterns within the data. As the layers narrow, heuristically, the network would refine these features into more sophisticated representations suitable for making predictions. The conjecture agrees with this common heuristic.

\begin{conjecture}[\cite{kileel2019expressive}, Conjecture 12]
    \label{conj:minimalWidthsUnimodal}
    Fix $L,d_0,d_L$ and $r$; any minimal (with respect to $\preccurlyeq$) vector of widths $\d=(d_0,d_1,\dots,d_L)$ such that the architecture $\mathcal{V}_{\d,r}$ is filling, is unimodal, i.e.\ there exists $i\in\{0,1,\dots,L\}$ such that $(d_0,\dots,d_i)$ is weakly increasing and $(d_{i},\dots,d_L)$ is weakly decreasing.
\end{conjecture}

\subsection{A look into the symmetries of an exceptional shallow network} \label{sec:symmetries-of-an-exceptional-shallow-network}

We explicitly describe a family of symmetries beyond multi-homogeneity (Lemma \ref{lem:multiHomog}) for the architecture $\mathbf{d}=(5,7,1),~r=3$, one of the exceptional cases of the Alexander--Hirschowitz Theorem. The expected dimension of $\V_{\mathbf{d},r}$ is 35, however, its actual dimension is 34. Therefore, there exists a one-dimensional family of symmetries in the parameter space not of the form as in Lemma \ref{lem:multiHomog}. \par 
In the following we will use an equivalent formulation of the Alexander--Hirschowitz Theorem as can be found in \cite[Theorem 1.1]{brambilla2008alexander}.

\begin{theorem}[Alexander--Hirschowitz, geometric version]
    Let $X$ be a general collection of $k$ double points in $\P^n$ and let $I_X(d)\subseteq \Sym_d(\C^n)$ be the subspace of polynomials through $X$, i.e., with all first partial derivatives vanishing at the points of $X$. Then the subspace $I_X(d)$ has the expected codimension $\min\{(n+1)k, \binom{n+d}{n}\}$ except in the cases as in Theorem \ref{thm:AH} (with the notation $d=r,~n=d_0-1,~k=d_1$).
\end{theorem}

Then geometrically the situation depicts as follows: consider seven points $p_1,\dots,p_7$ in $\P^4$. We are interested in cubics singular at these seven points, i.e.,\ we are looking for $f\in \Sym_3(\C^5)$ such that $f(p_i)=\mathrm{d}f_{p_i}=0$ for $i=1,\dots,7$. One would expect that no such cubic exists: the space of cubics in five variables has dimension 35 and the seven points are expected to impose 35 independent conditions. But through seven points there exists a rational normal curve $C$ which after a suitable coordinate transformation might be given as
\[
    C = \left\{(x_0:\dots :x_4)\in \P^4 \colon \rk \begin{pmatrix}
        x_0 & x_1 & x_2 \\
        x_1 & x_2 & x_3 \\
        x_2 & x_3 & x_4 \\
    \end{pmatrix} \leq 1\right\}.
\] 
A different way to write the curve $C$ is as follows. Assume the five points $p_1,\dots,p_5$ are the points $(1:0:0:0:0),\dots,(0:0:0:0:1)$; if $p_1,\dots,p_7$ are in general position $p_6$ and $p_7$ have non-zero coordinates then. Consider the birational Cremona transformation $\tau\colon (x_0:\dots:x_5) \mapsto (x_0^{-1} :\dots : x_5^{-1})$. Then $C$ is the preimage under $\tau$ of the line $\overline{\tau(p_6)\tau(p_7)}$.\par 
It is well-known that the secant variety of such a determinantal variety is given by
\[
    \sigma(C) = \left\{(x_0:\dots :x_4)\in \P^4 \colon \det \begin{pmatrix}
        x_0 & x_1 & x_2 \\
        x_1 & x_2 & x_3 \\
        x_2 & x_3 & x_4 \\
    \end{pmatrix} = 0\right\}.
\]
This is a cubic hypersurface with singular locus $C$; in particualar, $\sigma(C)$ is singular at $p_1,\dots,p_7$. See \cite[\S 3]{brambilla2008alexander}.\par
Now it is possible to vary $p_7$ along $C$ such that the construction of $\sigma(C)$ remains unchanged. Consider the PNN
\[
    W_2\circ \rho_3 \circ W_1 \mathbf{x}
\]
where $\mathbf{x}\in\R^5,~W_1\in\R^{7\times 5},~W_2\in\R^{1\times 7}$. Assume similar to above that the curve $\tau(C)$ is given by the $2\times2$ minors of 
\begin{equation}
    \label{equ:minorsDefC}
    \begin{pmatrix}
        w_{161}^{-1} & \dots & w_{165}^{-1} \\
        w_{171}^{-1} & \dots & w_{175}^{-1} \\
    \end{pmatrix}.
\end{equation}
We are looking for a symmetric matrix $G\in \mathrm{GL}(7,\R)$ such that $GW_1$ gives rise to the same curve $\tau(C)$. Again, w.l.o.g., we assume that the first five rows and columns of $G$ are the identity matrix. Thus, restricting to the last two rows and columns of $G$, we require that
\[
    \begin{pmatrix}
        g_1 & g_2 \\
        g_2 & g_3 \\
    \end{pmatrix}\begin{pmatrix}
        w_{161}^{-1} & \dots & w_{165}^{-1} \\
        w_{171}^{-1} & \dots & w_{175}^{-1} \\
    \end{pmatrix}
\]
have the same $2\times 2$ minors as \eqref{equ:minorsDefC}. Normalizing $G$ to have determinant one we can solve this to obtain
\[
    G = \begin{pmatrix}
        \mathrm{Id}_5 & 0 & 0 \\
        0 & g & \phi \\
        0 & \phi & \frac{1+\phi^2}{g} \\
    \end{pmatrix}
\]
where $\phi = (\frac{1}{2}(-1 + \sqrt{3}))^{1/2}$ and $g\in \R$ is a free parameter. Thus we have obtained a one-dimensional family of symmetries on the parameter space leaving the image of the parameter map $\Psi_{\mathbf{d},r}$ unchanged.

\section{Optimization} \label{sec:optimization}
In supervised machine learning tasks, a neural network is trained through the process of \emph{empirical risk minimization (ERM)}. Consider a training dataset $\mathcal{D} = \{ (\x_1, \y_1), (\x_2, \y_2), ..., (\x_N, \y_N) \}$, where each $\x_i$ is an input vector and $\y_i $ is the corresponding output or label. The objective in ERM is to find a function $f$ from a hypothesis space $\mathcal{H}$ that minimizes the empirical risk, which is defined as the average loss over the training data:
\[ \hat{R}(f) = \frac{1}{N} \sum_{i=1}^{N} {L}(\y_i, f(\x_i)). \]
Here, $\hat{R}(f)$ represents the empirical risk of the function $f$, ${L}$ is a loss function that measures the discrepancy between the predicted value $ f(\x_i) $ and the actual label $ \y_i $, and $N$ is the number of samples in the training dataset.\par

The goal of ERM is to select a function \( f^* \) from the hypothesis space \( \mathcal{H} \) that minimizes this empirical risk:
\begin{equation}
    f^* = \underset{f \in \mathcal{H}}{\arg\min} \ \hat{R}(f). \label{eq:minrisk}
\end{equation}

To minimize empirical risk, gradient-based optimization methods such as gradient descent (GD), stochastic gradient descent (SGD), and adaptive moment estimation (Adam) are commonly used. For convex objective functions, gradient-based algorithms can converge to the unique global minimum with appropriate hyperparameter choices. However, when optimizing neural networks, the loss landscape is typically non-convex. As a result, there is generally no theoretical guarantee that gradient-based methods will reach a global optimum $f^*$ for feedforward neural networks.

Recent research on convergence properties in gradient-based algorithms often focuses on the dynamics of the optimization process. One approach analyzes convergence rates by relaxing assumptions on the objective function. For example, it has been shown that for objective functions with a unique global minimum, gradient descent converges to the global minimum at a geometric rate if the function satisfies conditions such as the Polyak--Lojasiewicz inequality \cite{karimi2016linear, li2019convergence,khaled2020better,mei2021leveraging}, weak quasi-convexity \cite{hardt2018gradient}, or the restricted Secant inequality \cite{hardt2018gradient, yi2020exponential, guille2022gradient}. Another line of work explores algorithmic regularization. Empirical evidence suggests that large-scale neural networks used in practice are highly over-parameterized, often containing far more trainable parameters than training examples. Consequently, optimization objectives for such high-capacity models tend to have many global minima that interpolate the data perfectly. The choice of optimization algorithm thus introduces an implicit inductive bias, guiding the solution toward specific types of global minima. One of the first results in this direction shows that for almost all linearly separable datasets, gradient descent with any initialization and any bounded stepsize converge in direction to maximum margin separator with unit ${l}_2$ norm \cite{soudry2018implicit}. Both approaches build on the assumption of the existence and separation of global minima, while the global picture of the distribution of critical points is not well-understood. Through algebro-geometric approach, we study the geometry of the optimization landscape and provide the first explicit upper bound to the number of critical points when optimizing over a special family of polynomial neural networks.

Note that in the neural network setting, the hypothesis class $\mathcal{H}$ is just the corresponding neuromanifold, and the empirical minimization problem becomes 
\begin{equation}\label{eq:nnloss}
   {\ttheta}^* =\arg\min_{\ttheta}\ell_{F_{\ttheta}}\quad \text{where }\ell_{F_{\ttheta}}\coloneqq \frac{1}{N} \sum_{i=1}^{N} {\ell}(\y_i, F_{\ttheta}(\x_i))
\end{equation}
Considering the $\ell_2$ loss, we minimize the average distance between $F_{\ttheta}(\x)$ and $\y$.\ It turns out that, analogously to maximum likelihood estimation or Euclidean distance minimization, there exists a degree of this optimization problem constituting a measure for the algebraic complexity of solving this problem. In \S \ref{sec:learning-degree}, we study the static optimization properties of a PNN through the generic ED degree of the corresponding neurovariety.\ In \S \ref{sec:dynamic} and \S \ref{sec:backprop}, we explain the dynamic optimization process of gradient descent, followed by backpropagation.\ The algorithm biases induced by the optimization methods, known as \emph{implicit bias}, on general neural network architectures are still largely unknown.\ While limited theoretical results are known for the dynamic optimization process of neural networks in general, recent work has shown that gradient descent finds a global minimum of linear neural networks (i.e.,\ polynomial neural networks with $r=1$) under mild assumptions  \cite{arora2018convergence}.\ In \S \ref{sec:linear}, we review the trajectory-based convergence results on linear neural networks and phrase open questions for polynomial neural networks with activation degree $r\geq2$. Sections \S \ref{sec:dynamic}, \S \ref{sec:backprop} and \S \ref{sec:linear} are expository and do not contain new results.

\subsection{The Learning Degree of a PNN}\label{sec:learning-degree}

Let us consider the case of PNNs again. Often one tries to solve the optimization problem \eqref{eq:minrisk} using gradient based methods, see \S \ref{sec:dynamic}. Therefore, one is interested in computing the critical points on the neurovariety of the loss function $\ell_{p_{\w}}$ as defined in \eqref{eq:nnloss} for a PNN $p_{\w}$. In particular, the number of such points provides an upper bound for the number of functions the optimization process can converge to when using gradient descent.\par

\begin{definition}
    \label{def:learningDeg}
    If there is a finite number of critical points of $\ell_{p_{\w}}$ on the regular locus of $\V_{\d,r}$ which is constant for a general choice of training data $(\x_i,\y_i)_i$, we call this number the \emph{learning degree} of the network $p_{\w}$ with respect to loss $\ell$, and denote it by $\ldeg_{\ell}(p_{\w})$.
\end{definition}

In the following we are studying the case of the Euclidean loss function. We show that in this case the learning degree is independent of the choice of training data $(\x_i,\y_i)_i$. An upper bound of the learning degree is provided by the generic Euclidean distance (ED) degree. We would like to thank Matthew Trager for pointing out the connection to the generic ED degree.

\begin{definition}[{\cite[Definition 2.8]{MAG}}]
    Let $X$ be a variety, let $\Lambda = (\lambda_1,\dots,\lambda_n)\in\R^n_{+}$ and let $\mathbf{u}\in\R^n$ be a general point. The $\Lambda$-weighted Euclidean distance degree of $X$ is the number of (complex) critical points on $X$ of the function
    \[
        \|\x-\mathbf{u}\|_{\Lambda}^2 \coloneqq \sum_{i=1}^n \lambda_i(x_i-u_i)^2.
    \]

    For general weights $\Lambda$, this number is independent of $\Lambda$, and we call it the \emph{generic Euclidean distance (ED) degree} of $X$, which we denoted as $\mathrm{EDdeg}_{\mathrm{gen}}(X)$. 
\end{definition}

Let $(p_{\w}^{(i)}(\x))_i \in (\Sym_d(\R^n))^m$ be a tuple of polynomials in $n$ variables $\x = (x_1,\dots,x_n)$ occurring as the output of a fixed PNN with weight vector $\w$. Suppose we want the PNN to learn the tuple of polynomials $(f_i(\x))_i\in (\Sym_d(\R^n))^m$. To this end, we sample $N$ input vectors $\hat{\x}_1,\dots,\hat{\x}_N\in\R^n$ according to some probability distribution. We evaluate  $(f_i(\x))_i$ at these points to obtain $\hat{\y}_j \coloneqq (f_i(\hat{\x}_j))_i \in \R^{m}$ for $j=1,\dots,N$; the pairs $(\hat{\x}_j, \hat{\y}_j)_j$ constitute the training data.\par 
The learning process is defined by minimizing the average Euclidean distance between the $\hat{\y}_j$ and the points obtained from evaluating the network at the samples $\hat{\x}_j$. This process corresponds to the empirical risk minimization with $\ell_2$ loss.  Concretely, we have the following problem:
\begin{equation}
    \label{equ:training}
    \text{minimize } \frac{1}{N} \sum_{j=1}^N || (p_{\w}^{(i)}(\hat{\x}_j))_i - \hat{\y}_j ||^2 \text{ over weights } \w. 
\end{equation}
Note that $(p_{\w}^{(i)}(\hat{\x}_j))_i ,\hat{\y}_j\in\R^m$ for all $j=1,\ldots,N$, and we have
\[
    \frac{1}{N} \sum_{j=1}^N || (p_{\w}^{(i)}(\hat{\x}_j))_i - \hat{\y}_j ||^2=\frac{1}{N}\sum_{i=1}^m \sum_{j=1}^N || p_{\w}^{(i)}(\hat{\x}_j) - (\hat{\y}_j)_i ||^2. 
\]

Our goal is to show that for each $i=1,2,\ldots, m$,
this is a quadratic form in the coefficients of $p_{\w}^{(i)}$ and $f_i$ over training data.\par 

Write 
\[
    p_{\w}^{(i)}(\x) = \sum_{I\in \binom{[n]}{d}}\rho^i_I \x^I \quad\text{and}\quad f_i(\x) = \sum_{J\in\binom{[n]}{d}} \phi^i_{J}\x^J.
\]
We define a matrix $E$ as follows: $E$ is a block-diagonal matrix with blocks $E^i$ for $i=1,2,\dots,m$, and each block $E^i$ has as row and column labels the multiindices $\alpha,\beta\in\binom{[n]}{d}$. Then the entry $E^i_{\alpha,\beta}$ is defined as 
\[
    E^i_{\alpha,\beta} \coloneqq \frac{1}{N}\sum_{j=1}^N\hat{\x}_j^{\alpha+\beta}.
\]
Then we can write
\[\frac{1}{N} \sum_{j=1}^N || p_{\w}^{(i)}(\hat{\x}_j) - f_i(\hat{\x}_j) ||^2=\sum_{\alpha,\beta}(\rho^i_\alpha-\phi^i_\alpha)E^i_{\alpha,\beta}(\rho^i_\beta-\phi^i_\beta).\]
Therefore, the loss function is an $E^i_{\alpha,\beta}$-weighted Euclidean distance in the space $(\Sym_d(\R^n))^m$, the ambient space of the neurovariety $\V_{\d,r}$. For a general choice of training data $(\hat{\x}_i,\hat{\y}_i)_i$,
each block $E^i$ of the linear transformation is generic. Note, however, that all blocks $E^i$, for $i=1,\dots,m$, will be equal, hence, the transformation $E$ might be non-generic. We summarize the discussion in the theorem below.
\begin{theorem} \label{thm:learning-degree-is-upper-bounded-by-the-generic-ED-degree}
    The learning degree $\ldeg_{\ell_2}(p_{\w})$ of a PNN $p_{\w}$ with architecture $(\d, r)$ with respect to the Euclidean loss function exists. It is at most the generic Euclidean distance degree of its neurovariety, i.e.,
    \begin{equation}
        \label{equ:learningDegreeBound}
        \ldeg_{\ell_2}(p_{\w}) \leq \mathrm{EDdeg}_{\mathrm{gen}}(\V_{\d,r}).
    \end{equation}
    If $p_{\w}$ has only a single output neuron, one has equality in \eqref{equ:learningDegreeBound}.
\end{theorem}

Note that for large number of samples, the quadratic form defined by $E$ is in fact independent of the samples.

\begin{remark}
Suppose $\x_1,\ldots,\x_N$ are drawn independently from a fixed distribution $\mathcal{D}$ with known moments $\mu_k$, $k=1,2,\ldots$ By the law of large numbers, for any $\epsilon>0$,
\[\mathbb{P}(|E^i_{\alpha,\beta}-\mu_{\alpha+\beta}|>\epsilon)\rightarrow0\text{ as }N\rightarrow\infty.\] \par
Then for any $\epsilon>0$, 

\[\mathbb{P}\left(\left|\frac{1}{N} \sum_{j=1}^N || p_{\w}^{(i)}(\hat{\x}_j) - f_i(\hat{\x}_j) ||^2-(\boldsymbol{\rho}^i - \boldsymbol{\phi}^i)^T A(\boldsymbol{\rho}^i - \boldsymbol{\phi}^i)\right|>\epsilon\right)\rightarrow0\text{ as }N\rightarrow\infty.\] 
where $A$ is a fixed matrix of dimension $N\times N$ whose entries are determined by the distribution $\mathcal{D}$. 
\end{remark}

For the architecture $\mathbf{d}=(2,2,k),~r=2$, we compute the bound in \eqref{equ:learningDegreeBound} explicitly.

\begin{theorem}
\label{prop:EDded22k}
    For $k\geq 2$, the generic ED degree of $\V_{(2,2,k),2}$ is $8k^2 - 12k + 3$. Hence, the learning degree of the PNN with architecture $(\d=(2,2,k),r=2)$ is at most $8k^2 - 12k + 3$.
\end{theorem}

The generic ED degree of a variety can be computed as the sum of its polar degrees \cite[Corollary 2.14]{MAG}. For a smooth variety $X$, this can be reformulated in terms of degrees of Chern classes of $X$, see \cite[Theorem 4.20]{MAG}. However, $V \coloneqq \V_{(2,2,k),2}$ is not smooth: recall that $V$ consists of $k\times3$ matrices with rank $\leq 2$; this determinantal variety has a nonempty singular locus consisting of all matrices with rank $\leq 1$. Therefore, we have to make use of the more involved machinery of \emph{Chern--Mather classes}. The reader wishing to learn more about Chern classes in the context of optimization is referred to \cite[\S 4.3]{MAG}. We give a brief definition of Chern--Mather classes below, for more details the reader is referred to \cite{yokura1986polar}.\par 
Let $X\hookrightarrow \P^N$ be a projective variety of pure dimension $d$. Let $X_{\text{sm}}$ denote the open subvariety of $X$ of nonsingular points and consider the embedding
\[
    g\colon X_{\text{sm}} \hookrightarrow \Gr(d, T\P^N),\quad x\mapsto T_xX_{\text{sm}}.
\]
Then $\hat{X}\coloneqq \overline{\im(g(X_{\text{sm}}))}$ is called the \emph{Nash blow-up} of $X$. The restriction of the projection $\Gr(d, T\P^N)\rightarrow \P^N$ to $\hat{X}$ gives the \emph{Nash blow-up map} $\nu\colon \hat{X}\rightarrow X$. Restricting the tautological subbundle of $\Gr(d, T\P^N)$ to $\hat{X}$ gives the \emph{Nash tangent bundle} $\widehat{TX}$ of $\hat{X}$.

\begin{definition}
    The $i^{\text{th}}$ \emph{Chern--Mather class} $c^M_i(X)$ of $X$ is the pushforward
    \[
        c^M_i(X) := \nu_*\left(c_i(\widehat{TX} \cap [\hat{X}])\right).
    \]
    The total Chern--Mather class $c^M(X)$ is the sum of all $c^M_i(X)$.
\end{definition}

The following result expresses the generic ED degree of a determinantal variety in terms of degrees of Chern--Mather classes.

\begin{lemma}[{\cite[Proposition 5.5]{chernDetVar}}]
    \label{lem:EDdegCM}
    Let $V_{m,n,r}\hookrightarrow \P^{mn-1}$ be the determinantal variety of $m\times n$ matrices with rank $\leq r$. Then 
    \[
        \mathrm{EDdeg}_{\mathrm{gen}}(V_{m,n,r}) = \sum_{l=0}^{(m+k)(n-k)-1} \sum_{i=0}^{l} (-1)^i \binom{(m+k)(n-k)-i}{(m+k)(n-k)-l} \beta_{(m+k)(n-k)-1-i},
    \]
    where $k = m-r$ and $c^M(V_{m,n,r}) = \sum_{l=0}^{mn-1} \beta_l H^l \in \Z[H]/\langle H^{mn}\rangle = A(\P^{mn-1})$ with $H=c_1(\O_{\P^{mn-1}}(1))$.
\end{lemma}

The Chern--Mather class of a determinantal variety has a quite involved, yet very explicit description due to Zhang.

\begin{theorem}[{\cite[Theorem\ 4.3]{chernDetVar}}]\label{thm:22k}
    With the notation as in Lemma \ref{lem:EDdegCM},
    \[
        c^M(V_{m,n,r}) = \mathrm{trace}(A(m,n,r)\cdot \mathcal{H}(m,n,r)\cdot B(m,n,r)).
    \]
    Here, $A,B$ and $\mathcal{H}$ are the following $(m(n-k)+1) \times (m(n-k)+1)$ matrices (remember $k = m-r$):
    \begin{align*}
         & A(m,n,r)_{i,j} = \int_{\Gr(k,n)} c(\T_{\Gr(k,n)})c_i(\QQ^{\vee m})c_{j-i}(\SS^{\vee m}) \cap [\Gr(k,n)], \\
         & B(m,n,r)_{i,j} = \binom{m(n-k) - j}{i-j},\\
         & \mathcal{H}(m,n,r)_{i,j} = H^{mk+j-i},
    \end{align*}
    where $\SS$ and $\QQ$ are the universal sub- and quotient bundle of the Grassmannian, respectively, and $i,j=0,1,\dots,m(n-k)$.
\end{theorem}

\begin{proof}[Proof of Theorem \ref{prop:EDded22k}]
    We first need to compute the Chern--Mather class $c^M(V)=c^M(V_{k,3,2})$. Note that in this case $k=1$, so $\Gr(k,n)=\Gr(1,3)\cong\P^2$. Let $A(\P^2) = \Z[h]/\langle h^3\rangle$ be the Chow ring of $\P^2$ where $h=c_1(\O_{\P^2}(1))$. The universal subbundle $\SS$ is just $\O_{\P^2}(-1)$ and we have
    \[
        c(\T_{\P^2}) = 1+3h+3h^2,\quad c(\SS^{\vee k}) = 1 + kh + \frac{1}{2}k(k-1)h^2,\quad c(\QQ^{\vee k}) = 1-kh+\frac{1}{2}k(k+1)h^2.
    \]
    Using these, we compute, for example,
    \[
        A(k,3,2)_{1,2} = \int_{\P^2} (1+3h+3h^2)(-kh)kh \cap [\P^2] = -k^2.
    \]
    Similarly, we obtain that $A(k,3,2)$ is the $(2k+1)\times(2k+1)$ matrix
    \[
        A(k,3,2) = \begin{pmatrix}
            3 & 3k & \frac{1}{2}k(k-1) & 0 & \dots & 0 \\
            0 & -3k & -k^2 & 0 & \dots & 0 \\
            0 & 0 & \frac{1}{2}k(k+1) & 0 & \dots & 0 \\
            0 & 0 & 0 & 0 & \dots & 0 \\
            \vdots & \vdots & \vdots & \vdots & \ddots & \vdots \\
            0 & 0 & 0 & 0 & \dots & 0 \\
        \end{pmatrix}.
    \]
    The matrix $B(k,3,2)$ is given by
    \[
        B(k,3,2) = \begin{pmatrix}
            \binom{2k}{0} & 0 & 0 & \dots & 0 \\
            \binom{2k}{1} & \binom{2k-1}{0} & 0 & \dots & 0 \\
            \binom{2k}{2} & \binom{2k}{1} & \binom{2k}{0} & \dots & 0 \\
            \vdots & \vdots & \vdots & \ddots & \vdots \\
            \binom{2k}{2k} & \binom{2k-1}{2k-1} & \binom{2k-2}{2k-2} & \dots & \binom{0}{0} \\
        \end{pmatrix}
    \]
    and the matrix $\mathcal{H}(k,3,2)$ is
    \[
        \mathcal{H}(k,3,2) = \begin{pmatrix*}[l]
            H^k & H^{k+1} & H^{k+2} & \dots & 0 \\
            H^{k-1} & H^k & H^{k+1} & \dots & H^{3k-1} \\
            H^{k-1} & H^{k-1} & H^k & \dots & H^{3k-2} \\
            \vdots & \vdots & \vdots & \ddots & \vdots \\
            0 & 0 & 0 & \dots & H^k \\
        \end{pmatrix*},
    \]
    where $H=c_1(O_{\P^{3k-1}}(1))$ is the generator of $A(\P^{3k-1})$. Carrying out the multiplication $B(k,3,2)A(k,3,2)\mathcal{H}(k,3,2)$, one finds that the diagonal has entries
    \[
        \scalebox{0.9}{$\begin{pmatrix}
            \binom{2k}{0}\left( 3H^k + 3kH^{k-1} + \frac{1}{2}k(k-1)H^{k-2}\right) \\
            \binom{2k}{1}\left( 3H^{k+1} + 3kH^{k} + \frac{1}{2}k(k-1)H^{k-1}\right)+ \binom{2k-1}{0} \left(-3kH^k-k^2H^{k-1}\right) \\
            \binom{2k}{2}\left( 3H^{k+2} + 3kH^{k+1} + \frac{1}{2}k(k-1)H^{k}\right)+ \binom{2k-1}{1} \left(-3kH^{k+1}-k^2H^{k}\right) + \binom{2k}{0}\frac{1}{2}k(k+1)H^{k-1} \\
            \vdots \\
            \binom{2k}{2k}\cdot 0 + \binom{2k-1}{2k-1}\left(-3kH^{3k-1} - k^2H^{3k-2}\right) + \binom{2k-2}{2k-2}\frac{1}{2}k(k+1)H^{3k-3} \\
        \end{pmatrix}$}.
    \]
    Computing the trace we find that
    \begin{align*}
        c^M(V_{k,3,2}) & = \sum_{j=-2}^{2k-1} \left[
        3\binom{2k}{j} + 3k\left(\binom{2k}{j+1} - \binom{2k-1}{j}\right)+ \frac{1}{2}k(k-1)\binom{2k}{j+2}\right. \\
        & \left. + \frac{1}{2}k(k+1)\binom{2k-2}{j} - k^2\binom{2k-1}{j+1}
        \right]
    \end{align*}
    and hence, by Lemma \ref{lem:EDdegCM}, 
    \begin{multline*}
        \mathrm{EDdeg}_{\text{gen}}(\V_{(2,2,k),2}) = \sum_{l=0}^{2(k+1)-1} \sum_{i=0}^{l} (-1)^i\binom{2(k+1)-i}{2(k+1)-l} \left[ 3\binom{2k}{i-2}+ 3k\left(\binom{2k}{i-1} \right.\right.\\
        \left. \left.- \binom{2k-1}{i-2}\right) + \frac{1}{2}k(k-1)\binom{2k}{i} + \frac{1}{2}k(k+1)\binom{2k-2}{i-2} - k^2\binom{2k-1}{i-1}
        \right].
    \end{multline*}
    Rearranging the double sum to
    \begin{multline*}
        \sum_{i=0}^{2(k+1)-1} \left(\sum_{l=i}^{2(k+1)-1} \binom{2(k+1)-i}{2(k+1)-l}\right) (-1)^i\left[ 3\binom{2k}{i-2}+ 3k\left(\binom{2k}{i-1} \right.\right.\\
        \left. \left.- \binom{2k-1}{i-2}\right) + \frac{1}{2}k(k-1)\binom{2k}{i} + \frac{1}{2}k(k+1)\binom{2k-2}{i-2} - k^2\binom{2k-1}{i-1}
        \right]
    \end{multline*}
    and noting that
    \[
        \sum_{l=i}^{2(k+1)-1} \binom{2(k+1)-i}{2(k+1)-l} = 
        2^{2(k+1)-i} - 1,
    \]
    one verifies that the above expression indeed equals $8k^2 - 12k + 3$. 
\end{proof}

We would like to thank Mark Kong for providing the rearrangement argument of the double sum.

\begin{remark} \label{remark:learning-degree-223}
    For the architecture $(\d=(2,2,3),r=2)$, Theorem \ref{thm:22k} yields $\ldeg_{\ell_2}(p_{\w}) \leq 39$. A numerical computation for a random choice of block $E^i$ in the linear transformation $E$ suggests that the actual learning degree is $\ldeg_{\ell_2}(p_{\w}) \leq 3$. Moreover, we observe that all critical points of the loss function are actually real. This phenomenon should be further studied in future work.
\end{remark}

\subsection{Gradient-based optimization}\label{sec:dynamic}
In order to solve for \eqref{eq:minrisk}, one uses iterative optimization algorithms, including gradient descent and its variants, to find a local minimum of the objective function. In each step, the parameter vector $\ttheta$ is adjusted in the direction opposite to the gradient of the function at the current point, as the gradient points in the direction of the steepest ascent. The update rule can be expressed as  

\begin{equation}
    \ttheta^{(t+1)} = \ttheta^{(t)} - \eta \nabla_{\ttheta} \ell_{F_{\ttheta^{(t)}}}.
\end{equation}

The tuning parameter $\eta$ is referred to as the \emph{learning rate} of the optimization process. In practice, the value of $\eta$ is crucial as it determines how big a step is taken on each iteration. If $\eta$ is too small, the algorithm may take too long to converge; if it is too large, the algorithm may overshoot the minimum or fail to converge.\par 

For a polynomial neural network $p_\w$ with weights $\w=(W_1,W_2,\ldots,W_L)$, at time step $t$, each $W_i$ is updated by 
\begin{equation}
    W_i^{(t+1)} = W_i^{(t)} - \eta \nabla_{W_i} \ell_{p_{\w^{(t)}}}
\end{equation}
where $\nabla_{W_i} \ell_{p_{\w^{(t)}}}$ is computed efficiently by backpropagation. In subsection \ref{sec:backprop}, we explain the backpropagation algorithm, which can potentially find its applications in algebraic statistics beyond the neural network settings. 

\subsection{Experiments} 

We generated synthetic data consisting of $m = 5000$ datasets, each containing $N = 50$ data points. For each dataset, input pairs $\hat{\mathbf{x}} = (\hat{x}_1, \hat{x}_2)^\top$ were sampled uniformly from the range $[-1, 1]$. The corresponding outputs $\hat{\mathbf{y}} = (\hat{y}_1, \hat{y}_2, \hat{y}_3)^\top$ were generated using quadratic functions with randomly sampled coefficients:
\begin{align*}
    \hat{y}_1 &= c_{1,1} \hat{x}_1^2 + c_{1,2} \hat{x}_1 \hat{x}_2 + c_{1,3} \hat{x}_2^2, \\
    \hat{y}_2 &= c_{2,1} \hat{x}_1^2 + c_{2,2} \hat{x}_1 \hat{x}_2 + c_{2,3} \hat{x}_2^2, \addtocounter{equation}{1}\tag{\theequation}\label{equ:defSampleOutputs}\\
    \hat{y}_3 &= c_{3,1} \hat{x}_1^2 + c_{3,2} \hat{x}_1 \hat{x}_2 + c_{3,3} \hat{x}_2^2,
\end{align*} 
where the coefficients $c_{i,j}$ were sampled from a normal distribution $\mathcal{N}(0, 1)$.

We employ a polynomial neural network $p_{\w}$ with architecture $\d = (2,2,3)$ and quadratic activation function $(r=2)$. The parameters $\w$ consist of two weight matrices $\w = (W_1, W_2)$ with $W_1 \in \R^{2 \times 2}$ and $W_2\in \R^{3\times 2}$.

Each network was trained on its corresponding dataset using stochastic gradient descent (SGD) with an initial learning rate of $\eta = 0.1$. The mean squared error (MSE) loss function was minimized:
\begin{equation*}
    L(W_1, W_2) = \frac{1}{N} \sum_{i=1}^{N} \lVert p_{\w}(\hat{\x}^{(i)}) - \hat{\mathbf{y}}^{(i)} \rVert_2^2.
\end{equation*}
A learning rate scheduler reduced the learning rate by half every $1000$ epochs. Training proceeded for a maximum of $15000$ epochs or until the maximum gradient magnitude fell below a threshold of $1 \times 10^{-4}$, indicating convergence. Later, we refer to this hyperparameter as the gradient norm threshold. 

After training, we extracted the quadratic coefficients learned by each network to represent the functions it had approximated. The extraction was performed by analytically expressing the network's output as a quadratic function of the inputs. For each output unit $j$, the coefficients $(a_{1j}, a_{2j}, a_{3j})$ were computed as
\begin{align*}
    a_{1j} &= v_{j1} w_{11}^2 + v_{j2} w_{21}^2, \\
    a_{2j} &= 2 (v_{j1} w_{11} w_{12} + v_{j2} w_{21} w_{22}), \addtocounter{equation}{1}\tag{\theequation}\label{equ:coeff2fcts}\\
    a_{3j} &= v_{j1} w_{12}^2 + v_{j2} w_{22}^2,
\end{align*}
where $v_{jk}$ are entries of $W_2$ and $w_{kl}$ are entries of $W_1$.

To identify distinct functions learned by the network, we compared the extracted coefficients using a coefficient comparison method with a specified tolerance $\epsilon$. Two functions given by coefficients $a_{ij}^{(1)}$ and $a_{ij}^{(2)}$ as in \eqref{equ:coeff2fcts} were considered the same if
\begin{equation*}
    \left| a_{ij}^{(1)} - a_{ij}^{(2)} \right| < \epsilon \quad \forall i, j.
\end{equation*}

We systematically adjusted the tolerance $\epsilon$ to explore its impact on the number of distinct functions identified. By varying $\epsilon$, we could control the granularity of function distinction, allowing minor variations in coefficients to be considered functionally equivalent.

In addition, we performed perturbations on the function's coefficients to assess whether each learned function corresponded to a local minimum in the neuromanifold. Small perturbations $\delta$ were added to the coefficients,
\begin{equation*}
    \tilde{c}_{ij} = c_{ij} + \delta_{ij}, \quad \delta_{ij} \in [-\varepsilon, \varepsilon],
\end{equation*}
where $\varepsilon$ is a small value (e.g., $1 \times 10^{-4}$). The perturbed coefficients also give rise to a perturbed neural network which we denote by $\tilde{p}_{\w}$. The perturbed loss function then becomes
\begin{equation*}
    L_{\delta}(W_1, W_2) = \frac{1}{N} \sum_{i=1}^{N} \left\| \tilde{p}_{\w}(\hat{\x}^{(i)}) - \tilde{\mathbf{y}}^{(i)} \right\|_2^2,
\end{equation*}
where $\tilde{\mathbf{y}}^{(i)}$ are defined as in \eqref{equ:defSampleOutputs} but using the perturbed coefficients $\tilde{c}_{ij}$. If the value of the original loss function was less than or equal to the losses of all perturbed functions, we considered it to be a local minimum in the neuromanifold.

Through our experimental setup, we were able to
\begin{itemize}
    \item train $5000$ datasets with gradient norm with $15000$ epochs and gradient norm threshold $10^{-4}$;
    \item identify seventeen distinct functions whose coefficients differ by tolerance threshold $\epsilon=1/10$ and that were learned with frequency no less than $10$ out of $5,000$ datasets;
    \item check that for sixteen out of seventeen functions, the coefficient matrices of the functions have rank one; 
    \item verify the local minimality of the one learned function whose coefficient matrix has rank two in the neuromanifold by evaluating the loss landscape around the function.
\end{itemize}

\begin{figure}[H]
    \centering
    \includegraphics[width=\linewidth]{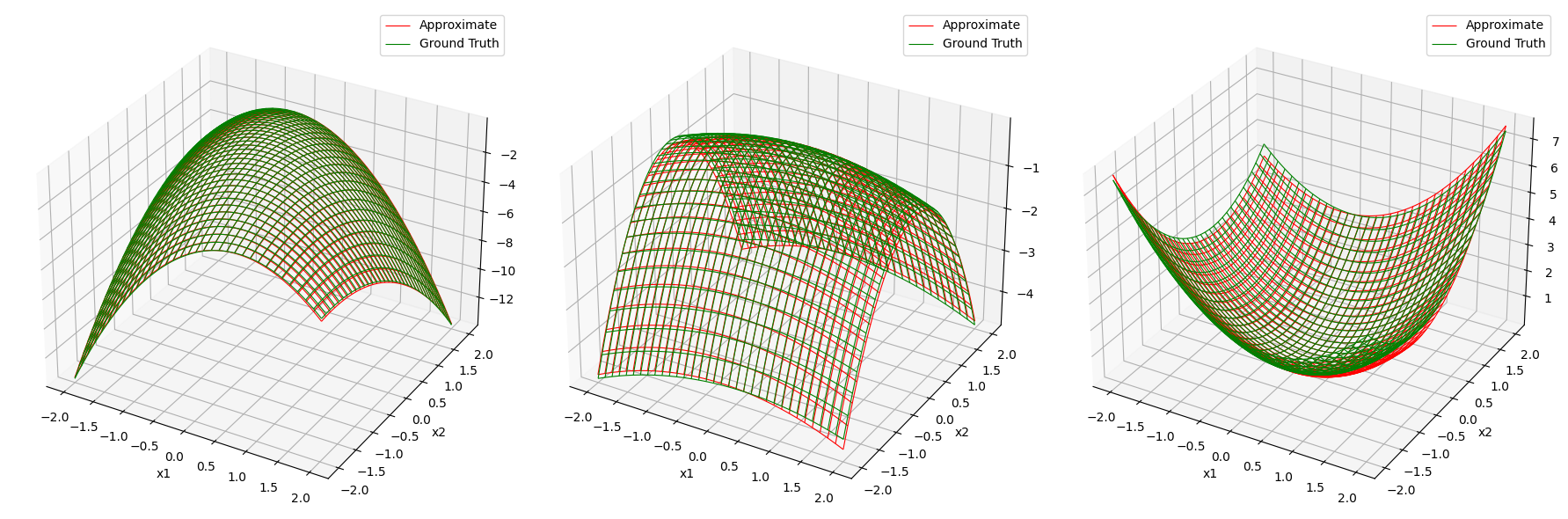}
    \caption{Comparison of the most frequent quadratic polynomial learned by a PNN with $\d=(2,2,3)$ and $r=2$ in each coordinate with the ground truth. The green manifold is generated based on the ground truth, while the red manifold is the most frequent function learned by the network.}
    \label{fig:function-comparison}
\end{figure}

The functions with rank-one coefficient matrices correspond to the singular points of the neuromanifold. Since the coefficients of quadratic functions in~\eqref{equ:defSampleOutputs} were generated randomly from the normal distribution, we expect the resulting function to be a nonsingular point of the neuromanifold. Therefore we consider only the learned function that is the nonsingular point of the neuromanifold. It is also the most frequently learned function by the network and it is depicted in Figure~\ref{fig:function-comparison}. Having one nonsingular local minimum is consistent with~\Cref{remark:learning-degree-223} that the learning degree is at most three.

The experiment implementation is available at 
\cite{mathrepo}. 
\begin{remark}
 The initialization of the network parameters poses challenges due to their tendency to suffer from vanishing and exploding gradients. As the depth of a polynomial neural network scales up, the degree of the polynomial increases and the network becomes sensitive to weight variations, causing gradients to either diminish or escalate exponentially during backpropagation. Poor initialization can lead to prolonged training times, convergence to suboptimal solutions, or even complete failure of the training process. While the general optimal initialization scheme remains an open research problem, in these experiments, we carefully initialize the weight and apply gradient clipping to regularize the training process.
 \end{remark}

\subsection{Backpropagation} \label{sec:backprop}
To compute the gradient of the loss function efficiently one uses the \emph{backpropagation algorithm}. Since Werbos introduced this algorithm in his 1974 PhD thesis \cite{werbos1974beyond}, backpropagation lies at the heart of modern successes in deep learning. For the computational algebraic geometer backpropagation is of interest as it provides a way to compute Jacobians of a parameter map $\Psi_{\d,r}$ exponentially faster than by na\"ive differentiation. We give a brief overview of the backpropagation algorithm following \cite[Chapter 2]{nielsen15NNDL} and explain how it can be used to compute Jacobian ranks.\par 
Let $\ell$ be a loss function (e.g.,\ as in \eqref{equ:training}) satisfying the following assumptions:
\begin{enumerate}
    \item $\ell$ is smooth;
    \item $\ell$ can be written as a sum of loss functions $\ell = \sum_{j=1}^N \ell_{\hat{\x}_j}$,
    where each $\ell_{\hat{\x}_j}$ depends only on the sample $\hat{\x}_j$;
    \item $\ell$ depends only on the output of the network $p_{\w}(\x)$, not on the state of intermediate layers; the training data $f(\hat{\x}_j)$ are considered as parameters.
\end{enumerate}
Backpropagation provides a way to compute the gradient $\nabla_{\w}\ell(\hat{\x}_1,\dots, \hat{\x}_N)$ efficiently. Note that in particular, if $\ell_{\hat{\x}}$ is just the network $F_{\ttheta}(\hat{\x})$ itself, this gives an algorithm for computing the gradient $\nabla_{\w}F_{\ttheta}(\hat{\x})$.\par 
Let $z^l_j$ denote the input into the $j^{\text{th}}$ neuron of the $l^{\text{th}}$ layer. The \emph{error} $\delta^l_j$ of this neuron is defined as 
\[
    \delta^l_j \coloneqq \frac{\del \ell}{\del z^l_j}.
\]

\begin{theorem}[{\cite{nielsen15NNDL}}]
    \label{thm:backprop}
    Let $a^l_j$ be the output of the $j^{\text{th}}$ neuron in the $l^{\text{th}}$ layer, and let $L$ denote the last layer of the network. Then the following equations hold:
    \[
        \delta^L_j = \frac{\del \ell}{\del a^L_j}\sigma'_{L,j}(z_j^L), \quad 
        \delta^l_j = \sigma'_{l,j}(z^l_j)\sum_{k=1}^{d_{l}} w_{l+1,j,k} \delta^{l+1}_{k},\quad
        \frac{\del \ell}{\del w_{l,j,k}} = a^{l-1}_k \delta^l_j.
    \]
\end{theorem}

\begin{proof}
    These are simple consequences of applying the chain rule.
\end{proof}

These equations give rise to Algorithm \ref{alg:backpropagation} below. The idea is to first compute the outputs of each neuron $a^l_j$ via a ``forward'' run through the network and then compute the errors $\delta^l_j$ using the first two equations of Theorem \ref{thm:backprop} running ``backwards'' through the network. Finally, the gradient can be computed using the third equation from Theorem \ref{thm:backprop}.

\begin{algorithm}
    \caption{Backpropagation}
    \label{alg:backpropagation}
    \hspace*{-3em} \textbf{Input:} A neural network $F_{\ttheta}$, a cost function $C$ with single input sample $\hat{\x}$ \\
    \hspace*{-14em} satisfying the assumptions above \\
    \hspace*{-22.2em} \textbf{Output:} The gradient $\nabla_{\w}\ell(\hat{\x})$\\
    \begin{algorithmic}[1]
    \State $a^0 \gets \hat{\x}$
    \For{$l$ from $1$ to $L$:}
        \Comment{forward loop}
        \State $z^l \gets W_l a^{l-1}$
        \State $a^l \gets \sigma_l(z^l)$
    \EndFor
    \State $\delta^L \gets \nabla_{\mathbf{a}}\ell \odot \sigma'_L(z^L)$
    \Comment{$\odot$ denotes the Hadamard (elementwise) product}
    \For{$l$ from $L-1$ to 1:}
        \Comment{backward loop}
        \State $\delta^l \gets (W_l^T \delta^{l+1}) \odot \sigma'_l(z^l)$
    \EndFor
    \State \Return $\nabla_{\w}\ell(\x) = (a^{l-1}_k \delta^l_j)_{j,k,l}$
    \end{algorithmic}
\end{algorithm}

If we choose the loss function $\ell$ to be the network $F_{\ttheta}$ itself, backpropagation returns
the gradient $\nabla_{\w}F_{\ttheta}(\hat{\x})$. However, when we are interested in computing the Jacobian of $\Psi_{\d,r}$, we need to compute the partial derivatives
$\del_{w_{j,k,l}} c^{(i)}_I$,
whereas
\begin{equation}
    \label{equ:backpropLinSys}
    \frac{\del F_{\ttheta}^{(i)}}{\del w_{j,k,l}}(\hat{\x}) = 
    \sum_{I} \frac{\del c^{(i)}_I}{\del w_{j,k,l}} \hat{\x}^I.
\end{equation}
By choosing $N\coloneqq \binom{r^{L-1}+d_0-1}{d_0-1}$ many samples $\hat{\x}_1,\hat{\x}_2,\dots, \hat{\x}_N$ and applying backpropagation to each of them, we obtain a linear system from \eqref{equ:backpropLinSys} in the unknowns $\frac{\del c^{(i)}_I}{\del w_{j,k,l}}$. For a generic choice of samples, this system has a unique solution and we obtain the Jacobian $J_{\Psi_{\d,r}}$. \par 
This routine allows us to compute Jacobians and in particular the dimension of the neurovariety exponentially faster than via the na\"ive method of computing all derivatives of the parametrization $\Psi_{\d,r}$. It can be easily adapted for computations over finite fields which makes dimension computations even faster. In \cite{kileel2019expressive} this routine has already been implemented, however (as of Oct 9 2024), their implementation at \cite{KTBgithub} contains a mistake in the set-up of the linear system described above. We provide a corrected implementation at \cite{mathrepo}.

\subsection{Case study: Linear Neural Networks} \label{sec:linear}
\vspace{0.4em} 

In addition to static properties of neural network optimization, we are also interested in the trajectory of practical optimization algorithms. This problem presents a more complex challenge, particularly in its general form. Limited results are known for the special subfamily of polynomial neural networks with activation degree $r=1$. In this subsection, we survey prior results on the trajectory of the gradient descent algorithms applied to linear neural networks from literature in the algebraic statistics and machine learning communities separately. An interesting open question is to extend the analysis to polynomial neural networks with general activation degrees.

Linear neural networks are a special class of polynomial neural networks with the activation function  $\rho_1(\x)=\x$. The associated map from the set of parameters to vectors of homogeneous polynomials of degree one in $d_0$ many variables is given by
\[
\Psi_{\d,1}: \R^{d_0 \times d_1} \times \cdots \times \R^{d_{L-1} \times d_{L}} \rightarrow (\Sym_{1}(\R^{d_0}))^{d_L}, \quad \w \mapsto p_{\w}=
\begin{bmatrix}
p_{\w}^{(1)}\\
\vdots\\
p_{\w}^{(d_L)}
\end{bmatrix}.
\]
The neuromanifold for the architecture $\d$ is
\[\M_{\d,1}=\{M\in\mathbb{R}^{d_0\times d_L}\mid \rk(M)\leq\min\{d_0,d_1,\ldots,d_L\}\}.\]
Note that in the linear case the neuromanifold is always equal to the neurovariety.

We denote the map from $\R^{d_0\times \cdots\times d_L}$ to $\R^{d_0\times d_L}$ as 
\[\mu_{\d}(W_1,W_2,\ldots,W_L)=W_L\cdots W_2W_1.\]
\begin{proposition}[{\cite[Theorem\ 4]{trager2019pure}}]
Let $\tilde{d}=\min(\d)$, $\w=(W_1,\ldots,W_L)$ and $\overline{W}=\mu_{\d}(\w)$.
\begin{itemize}
    \item(Filling case) If $\tilde{d}=\min\{d_0,d_L\}$, the differential $d\mu_d(\w)$ has maximal rank equal to $\dim (\M_{\d,1})=d_0 d_L$ if and only if for any $i\in[L-1]$,  either $\rk(W_j)=d_L$ for all $j>i$ or  $\rk(W_j)=d_0$ for all $j<i+1$. 
    \item(Non-filling case) If $\tilde{d}<\min\{d_0,d_L\}$, the differential $d\mu_d(\w)$ has maximal rank equal to $\dim (\M_{\d,1})=\tilde{d}(d_0+d_L-\tilde{d})$ if and only if $\rk(\overline{W})=\tilde{d}$. 
\end{itemize}
\end{proposition}

Besides static properties of the loss landscape, in the realm of neural networks, optimization often focuses on devising and refining algorithms to efficiently find the best model parameters that minimize a loss function.\ Key aspects of optimization include developing and analyzing optimization algorithms like gradient descent and its variants, understanding the impact of different loss functions, and ensuring model generalization through regularization.\  Additionally, researchers delve into convergence analysis to assess how and when algorithms reach optimal solutions. Since neural networks are, in general, complex models, the loss landscape is often non-convex, resulting in a significant risk of converging to a suboptimal local minimum. Analyzing the convergence trajectory is one of the most challenging problems in machine learning theory. We recall results on the convergence trajectory of deep linear neural networks with respect to $\ell_2$ loss in \Cref{thm:linearconv}.\ This analysis provides crucial insights into the efficiency and effectiveness of how learning algorithms perform on these architectures.

\begin{theorem}[{\cite[Theorem\ 1]{arora2018convergence}}, informal]\label{thm:linearconv}
Consider training deep linear neural networks with $\ell_2$ loss using the gradient descent algorithm.\ Assume that at initialization, the weight matrices are well-conditioned.\ In addition, there exists $\delta>0$ such that at any time step $T$, 
\begin{equation}
    \|W_{j+1}^T W_{j+1} - W_j W_j^T\|_F \leq \delta, 
\end{equation}
for all $j=1,2,\ldots,L-1$. For any $\epsilon>0$, with a sufficiently small learning rate, there exists $T_0$ such that the loss is no greater than $\epsilon$ for any $T\geq T_0$.  
\end{theorem}

Theorem \ref{thm:linearconv} shows that, when minimizing the $\ell_2$ loss for a deep linear network over a dataset with zero mean and unit variance, gradient descent converges to the global minimum at a linear rate, given mild assumptions. Similar results have been shown for deep linear neural networks with losses beyond the $\ell_2$ loss. In \cite{brechet2023critical}, the authors consider deep linear neural networks trained with the Bures--Wasserstein distance, a loss function commonly used in generative models. They analyze not only the critical points of the loss landscape but also the convergence of both gradient flow and gradient descent for the Bures--Wasserstein loss.

The convergence trajectory of polynomial neural networks with activation degree $r \geq 2$ remains an open problem.\ Heuristically, having a higher degree of activation reduces the symmetry in the parameter space, which may lead to interesting geometric properties in the trajectory of various optimization algorithms. We aim to establish theoretical guarantees for the convergence of general polynomial neural networks in our future work.

\section{Conclusion}
We have studied the functional space of neural networks with polynomial activation functions from an algebro-geometric perspective. These networks, characterized by their ability to model complex, high-degree polynomial relationships in data, provide a unique lens through which the intricate dynamics of deep learning models can be better understood. Moreover, since a large class of functions can be efficiently approximated by polynomials, delving into the polynomial framework 
can be important for the development of more efficient and accurate models.

In this paper, we have focused on the \emph{geometry} of polynomial neural networks. We have characterized neuromanifolds  and neurovarieties for some specific architectures, and studied their dimensions.
The geometry and dimension of the neuromanifold relates directly to the expressive power of the neuromanifold. Beyond expressivity, we also studied the optimization process of polynomial neural networks, and bound
the number of critical points of a loss function using the \emph{learning degree}. Besides our conjectures from \S \ref{sec:dimension} concerning the dimension of neurovarieties, there are further questions to be studied, for example: What are the learning degrees of more general architectures? How many of the critical points of the loss functions are real?

\bibliographystyle{alpha}
\bibliography{bibliography}

\end{document}